\documentclass[11pt,a4paper]{amsart}
\usepackage[a4paper,margin=1in]{geometry}
\usepackage[T1]{fontenc}
\usepackage{lmodern}
\usepackage[utf8]{inputenc}
\usepackage{amsmath,amssymb,amsthm,mathtools}
\usepackage{mathrsfs}
\usepackage{enumitem}
\usepackage[protrusion=true,expansion=false]{microtype}
\usepackage{graphicx}
\usepackage{float}
\usepackage{placeins}
\usepackage{booktabs}
\usepackage[hidelinks]{hyperref}

\theoremstyle{plain}
\newtheorem{theorem}{Theorem}[section]
\newtheorem{lemma}[theorem]{Lemma}
\newtheorem{proposition}[theorem]{Proposition}

\newtheorem{assumption}[theorem]{Assumption}
\theoremstyle{definition}
\newtheorem{definition}[theorem]{Definition}
\theoremstyle{remark}
\newtheorem{remark}[theorem]{Remark}

\newcommand{\Ga}{\Gamma}
\newcommand{\D}{\mathcal D}
\newcommand{\A}{\mathcal A}
\newcommand{\C}{\mathbb C}

\newcommand{\N}{\mathbb N}
\newcommand{\PV}{\operatorname{p.v.}}
\newcommand{\dist}{\operatorname{dist}}

\newcommand{\spanop}{\operatorname{span}}

\newcommand{\Log}{\operatorname{Log}}
\newcommand{\diag}{\operatorname{diag}}

\newcommand{\eps}{\varepsilon}

\title[Regularized B-spline--Heaviside collocation]{A Regularized B-Spline--Heaviside Collocation Method for Cauchy Singular Integral Equations with Piecewise H\"older Solutions}

\author{Maria Capcelea}
\address{Institute of Mathematics and Computer Science, Moldova State University, Chisinau, Republic of Moldova}
\email{maria.capcelea@usm.md}

\author{Titu Capcelea}
\address{Department of Computer Science, Moldova State University, Chisinau, Republic of Moldova}
\email{titu.capcelea@usm.md}

\date{}

\begin{document}
\raggedbottom
\setlength{\intextsep}{9pt plus 2pt minus 2pt}
\setlength{\textfloatsep}{10pt plus 2pt minus 2pt}
\setlength{\floatsep}{9pt plus 2pt minus 2pt}

\begin{abstract}
We develop a B-spline--Heaviside collocation method for Cauchy singular integral equations on a smooth closed $C^2$ contour when the exact solution is piecewise H\"older continuous with finitely many prescribed jumps.  Since the Cauchy singular integral of a discontinuous function generally has logarithmic terms at the jump points, we study
\[
        M=cI+dS+K:X_\alpha\to Y_\alpha,
        \qquad
        X_\alpha=PH^\alpha(\Ga,\D),
        \quad
        Y_\alpha=PH^\alpha_{\log,*}(\Ga,\D),
\]
where $Y_\alpha$ is a logarithmically enlarged piecewise H\"older space with lateral H\"older logarithmic coefficients.  The discontinuous component is represented by a nonredundant system of normalized relative Heaviside functions adapted to the closed contour. Collocation uses point evaluations at spline nodes separated from the jump set together with logarithmic-coefficient functionals at the jumps.  Assuming continuous stability of $M:X_\beta\to Y_\beta$, mesh-uniform scaled discrete stability of the regularized collocation operators, and a scaled consistency estimate for exact-jump approximants, we prove existence and uniqueness for sufficiently fine meshes and the error bound
\[
        \|\varphi-\varphi^H_{n_B}\|_{X_\beta}
        \le C h_B^{\alpha-\beta}\|\varphi\|_{X_\alpha},
        \qquad 0<\beta<\alpha<1.
\]
We give a matrix realization of the regularized scheme, including the logarithmic and point-collocation blocks, singularity-subtracted evaluation of the Cauchy action on splines, principal-value-safe arc formulas for the Heaviside terms, and an implementation algorithm. An abstract perturbation result shows that the same rate is preserved under sufficiently accurate quadrature. Numerical experiments with arcwise errors and finite-dimensional stability and consistency indicators support the theoretical assumptions over the tested meshes.
\end{abstract}

\maketitle

\section{Introduction}

Cauchy singular integral equations on closed contours form a traditional class of equations in complex analysis, boundary integral methods, and mathematical physics.  They arise in boundary formulations of elliptic problems, in interface and transmission models, in crack problems, in plane elasticity, in potential theory, and in many models where a boundary condition changes from one part of the contour to another.  For equations with smooth coefficients and smooth right-hand sides, polynomial, trigonometric, spline, Nystr\"om, and Galerkin discretizations are well understood and can provide stable and high-order approximations; see, for example, the classical and numerical references~\cite{Muskhelishvili,Gakhov,ProssdorfSilbermann,Atkinson,JunghannsSilbermann}.

The numerical situation changes substantially when the coefficients, the data, or the solution are only piecewise regular.  In many applications the boundary is divided into finitely many arcs, and the physical model prescribes different coefficients or different boundary conditions on different arcs.  The resulting unknown may be H\"older continuous on each arc but may have finite jumps at the interface points.  Spline collocation and Galerkin methods for singular integral equations with piecewise continuous coefficients have been studied in related settings in~\cite{ProssdorfRathsfeld1984,ProssdorfRathsfeld1986}.  Smooth global approximation spaces are not structurally adapted to such functions.  If a discontinuity is approximated indirectly by polynomials, trigonometric polynomials, or globally continuous splines, spurious oscillations and a loss of strong convergence near the jump points are typical.  A natural remedy is to build the discontinuity structure into the trial space itself.

The approximation space used in this paper combines periodic B-splines with finitely many step functions.  The B-spline part approximates the continuous H\"older component of the solution, while the Heaviside part represents the jumps explicitly.  This idea is particularly useful when the jump set is known in advance.

The assumption that the jump set is prescribed is natural in many models: material interfaces, switching boundary conditions, crack tips, junctions, and other geometric or physical features usually determine the discontinuity locations before discretization begins.  If the jump locations are not known a priori, they may be estimated from discrete information about the values of the function on the contour and then used as input for the enriched collocation scheme.  Related techniques for recovering distinguished points of functions from samples on simple closed contours were developed, for instance, in ~\cite{CapceleaCapcelea2020}.

There is, however, an important analytical feature specific to Cauchy singular integral equations.  If $\varphi$ is piecewise H\"older and has a nonzero jump at a point $t_j^d$, then the Cauchy singular integral $S\varphi$ generally contains a logarithmic singularity at $t_j^d$.  Consequently, in the natural strong piecewise H\"older norm with finite one-sided limits, one does not in general have $S:PH^\alpha(\Ga,\D)\to PH^\alpha(\Ga,\D)$.
Thus an equation of the form
\begin{equation*}
        c(t)\varphi(t)+d(t)(S\varphi)(t)+(K\varphi)(t)=f(t),\qquad t\in\Ga\setminus\D,
\end{equation*}
should not be interpreted, in the discontinuous case, as an equation in a single piecewise H\"older space.  The image of $S\varphi$ belongs naturally to a logarithmically enlarged space.

For this reason we formulate the analysis in a pair of spaces
\[
        X_\alpha:=PH^\alpha(\Ga,\D),\qquad
        Y_\alpha:=PH^\alpha_{\log,*}(\Ga,\D).
\]
The space $Y_\alpha$ contains functions that are H\"older on each continuity arc after subtracting endpoint logarithms.  The logarithmic coefficients are allowed to be lateral H\"older functions, not merely constants.  This is essential because the singular term appears multiplied by the coefficient $d\in X_\alpha$.  The product of a H\"older coefficient and a logarithmic singularity has a H\"older coefficient multiplying the logarithm.  Hence the crucial structural property is $X_\alpha\cdot Y_\alpha\subset Y_\alpha$.
The space $Y_\alpha$ is not required to be an algebra; indeed, the product of two logarithmic functions would generally produce squared logarithms.  The module property above is exactly what is needed for the operator $cI+dS+K$.

Another structural point concerns the use of Heaviside functions on a closed contour.  A single-valued step function on a closed curve cannot have exactly one nonzero jump, because the sum of jumps after one complete circuit must be zero.  We therefore use relative Heaviside functions.  A base jump point $t_1^d$ is fixed, and the relative Heaviside function $G_j$, $j=2,\ldots,n_d$, has jump $+1$ at $t_j^d$ and jump $-1$ at $t_1^d$.  Thus the discontinuous enrichment has dimension $n_d-1$, and the full trial space has dimension $n_B+n_d-1$.

The main convergence theorem is conditional in the standard numerical-analysis sense.  We assume a continuous stability estimate
\[
        \|u\|_{X_\beta}\le C\|Mu\|_{Y_\beta},
        \qquad 0<\beta<\alpha<1,
\]
and a corresponding scaled uniform discrete stability estimate for the eliminated collocation functionals.  The scaling is natural in the $X_\beta$-norm because point residuals control a H\"older seminorm only after the standard inverse estimate for splines is taken into account.  Under these assumptions and the associated scaled consistency estimate, the B-spline--Heaviside collocation solution converges to the exact solution in $X_\beta$ with rate $O(h_B^{\alpha-\beta})$.

The paper is organized as follows.  Section~\ref{sec:spaces} introduces the contour, the piecewise H\"older spaces $X_\alpha$, and the logarithmic spaces $Y_\alpha$.  Section~\ref{sec:operator} proves the mapping properties of the Cauchy and regular integral operators and discusses when a logarithmic right-hand side can be avoided.  Section~\ref{sec:relative} develops the relative Heaviside decomposition on a closed contour, the B-spline--Heaviside approximation estimate, and the matrix-level computational realization of the regularized method.  Section~\ref{sec:collocation} formalizes the regularized collocation scheme in terms of bounded observation functionals.  Section~\ref{sec:exactconv} proves a conditional existence, uniqueness, and convergence theorem for the exact operator action.  Section~\ref{sec:quadrature} proves an abstract perturbation theorem for quadrature-based implementations.  Section~\ref{sec:numerical} presents numerical experiments, including untrimmed maximum and discrete piecewise H\"older-type error diagnostics, together with inverse-norm and scaled-consistency indicators associated with the discrete hypotheses.  Section~\ref{sec:conclusion} contains the conclusions.

Throughout the paper, $C$ denotes a positive constant that may change from line to line but is independent of the mesh size.

\section{Piecewise H\"older and logarithmic spaces}\label{sec:spaces}

\subsection{Geometry of the contour}

Let $\Ga\subset\C$ be a simple, closed, positively oriented contour of class $C^2$.  We assume that $\Ga$ admits a $2\pi$-periodic parametrization
\[
        \gamma:[0,2\pi]\to\Ga,
        \qquad \gamma(0)=\gamma(2\pi),
\]
with $\gamma\in C^2$ and $\gamma'(\theta)\ne0$ for all $\theta$.  The $C^2$ regularity ensures that the geometric remainder in the local logarithmic expansion of the Cauchy transform is $C^1$; consequently it belongs to $H^\alpha$ for every $0<\alpha<1$.  The Euclidean distance on $\Ga$, the arclength distance, and the parameter distance are locally equivalent.  Therefore H\"older regularity on $\Ga$ may be described equivalently in the contour variable or in the parameter variable.  Throughout the paper the exponents satisfy
\[
        0<\beta<\alpha<1.
\]

Let
$
        \D=\{t_1^d,\ldots,t_{n_d}^d\}\subset\Ga
$
be a finite set of distinct points, ordered according to the orientation of $\Ga$.  These points split $\Ga$ into the open arcs
$
        \A(\Ga\setminus\D)=\{\Ga_1,\ldots,\Ga_{n_d}\},
$
where $\Ga_k$ is the oriented open arc from $t_k^d$ to $t_{k+1}^d$ and $t_{n_d+1}^d=t_1^d$.  We call $\D$ the jump set.

For a function $u$ with one-sided limits at a jump point $t_j^d$, we write
\[
        [u]_{t_j^d}:=u(t_j^d+0)-u(t_j^d-0),
\]
where the signs are understood with respect to the positive orientation.  We use the left-continuous representative, so that $u(t_j^d)=u(t_j^d-0)$ when a point value is needed.

\subsection{The space $X_\alpha=PH^\alpha(\Ga,\D)$}

\begin{definition}\label{def:Xalpha}
For $0<\alpha<1$, the space $X_\alpha=PH^\alpha(\Ga,\D)$ consists of all functions $u:\Ga\to\C$ such that:
\begin{enumerate}[label=(\roman*)]
\item $u$ is H\"older continuous of exponent $\alpha$ on each arc $\Ga_k$;
\item the one-sided limits $u(t_j^d\pm0)$ exist and are finite for all $t_j^d\in\D$;
\item $u$ is represented by its left-continuous values at the points of $\D$.
\end{enumerate}
The norm is
\begin{equation*}
        \|u\|_{X_\alpha}:=
        \max_{1\le k\le n_d}
        \|u\|_{H^\alpha(\Ga_k)},
        \qquad
        \|u\|_{H^\alpha(\Ga_k)}:=\|u\|_{\infty,\Ga_k}+[u]_{\alpha,\Ga_k},
\end{equation*}
where
\[
        \|u\|_{\infty,\Ga_k}:=\sup_{t\in\Ga_k}|u(t)|,
   \qquad
        [u]_{\alpha,\Ga_k}:=
        \sup_{\substack{t,s\in\Ga_k\\ t\ne s}}
        \frac{|u(t)-u(s)|}{|t-s|^\alpha}.
\]
If the supremum is taken over the corresponding closed arc, the endpoint values are understood as the appropriate one-sided limits.
\end{definition}

The space $X_\alpha$ is a Banach space.  If $0<\beta<\alpha<1$, then $X_\alpha$ is continuously embedded into $X_\beta$.

\begin{lemma}[Algebra property]\label{lem:Xalgebra}
For $0<\alpha<1$, $X_\alpha$ is a Banach algebra.  More precisely,
\[
        \|uv\|_{X_\alpha}\le C\|u\|_{X_\alpha}\|v\|_{X_\alpha},
        \qquad u,v\in X_\alpha.
\]
\end{lemma}

\begin{proof}
On every arc $\Ga_k$, the classical H\"older product estimate gives
\[
        [uv]_{\alpha,\Ga_k}
        \le \|u\|_{\infty,\Ga_k}[v]_{\alpha,\Ga_k}
        +\|v\|_{\infty,\Ga_k}[u]_{\alpha,\Ga_k}.
\]
Together with the corresponding sup-norm estimate and the maximum over finitely many arcs, this proves the assertion.
\end{proof}

\subsection{The logarithmic target space $Y_\alpha$}

The space $Y_\alpha$ is designed to contain logarithmic endpoint singularities on the arcs of $\Ga\setminus\D$.  Let $s$ be arclength on $\Ga_k$, $0<s<\ell_k$, so that $s=0$ corresponds to $t_k^d$ and $s=\ell_k$ corresponds to $t_{k+1}^d$.  Choose once and for all smooth cut-off functions $\chi_k^-$ and $\chi_k^+$ on $[0,\ell_k]$ such that $\chi_k^-=1$ near $0$ and vanishes away from $0$, while $\chi_k^+=1$ near $\ell_k$ and vanishes away from $\ell_k$.

\begin{definition}[Logarithmic piecewise H\"older space]\label{def:Yalpha}
The space $Y_\alpha=PH^\alpha_{\log,*}(\Ga,\D)$ consists of all functions $u$ on $\Ga\setminus\D$ such that on every arc $\Ga_k$ one can write
\begin{equation}\label{eq:Yrepresentation}
        u(t_k(s))=u_{k,0}(s)
        +\chi_k^-(s)u_{k,-}(s)\log s
        +\chi_k^+(s)u_{k,+}(s)\log(\ell_k-s),
\end{equation}
where
\[
        u_{k,0},u_{k,-},u_{k,+}\in H^\alpha[0,\ell_k].
\]
The norm is defined by the infimum over all representations of the form \eqref{eq:Yrepresentation}:
\begin{equation*}
        \|u\|_{Y_\alpha}:=
        \max_{1\le k\le n_d}
        \inf\left(
        \|u_{k,0}\|_{H^\alpha[0,\ell_k]}
        +\|u_{k,-}\|_{H^\alpha[0,\ell_k]}
        +\|u_{k,+}\|_{H^\alpha[0,\ell_k]}
        \right).
\end{equation*}
\end{definition}

The star in $PH^\alpha_{\log,*}$ emphasizes that the logarithmic coefficients are lateral H\"older functions.  This is stronger and more flexible than allowing only constant logarithmic coefficients.

\begin{lemma}\label{lem:Ybanach}
The space $Y_\alpha$ is a Banach space.  Moreover, $X_\alpha\hookrightarrow Y_\alpha$ continuously.
\end{lemma}

\begin{proof}
The representation \eqref{eq:Yrepresentation} identifies $Y_\alpha$ on each arc with a quotient of the product space
\[
        H^\alpha[0,\ell_k]\times H^\alpha[0,\ell_k]\times H^\alpha[0,\ell_k]
\]
by the kernel of the linear map that sends the triple of coefficients to the represented function.  Since the product space is Banach and the norm is the quotient norm, the arcwise space is Banach.  A finite product of these arcwise spaces is Banach, which proves completeness of $Y_\alpha$.  The embedding $X_\alpha\hookrightarrow Y_\alpha$ is obtained by choosing the logarithmic coefficients equal to zero.
\end{proof}

Although the coefficient functions in \eqref{eq:Yrepresentation} need not be unique away from the endpoints, their endpoint values corresponding to the leading logarithmic terms are intrinsic.  We denote the outgoing and incoming logarithmic coefficient functionals by
\begin{align*}
        \lambda_j^+(u)&:=\text{coefficient of }\log|t-t_j^d|\text{ on the arc leaving }t_j^d,\\
        \lambda_j^-(u)&:=\text{coefficient of }\log|t-t_j^d|\text{ on the arc entering }t_j^d.
\end{align*}
These are bounded linear functionals on $Y_\alpha$.

\begin{lemma}[Module property]\label{lem:moduleproperty}
Let $0<\alpha<1$.  If $a\in X_\alpha$ and $u\in Y_\alpha$, then $au\in Y_\alpha$ and
\[
        \|au\|_{Y_\alpha}
        \le C\|a\|_{X_\alpha}\|u\|_{Y_\alpha}.
\]
\end{lemma}

\begin{proof}
On an arc $\Ga_k$, write $u$ in the form \eqref{eq:Yrepresentation}.  Since $a|_{\Ga_k}\in H^\alpha[0,\ell_k]$ and $H^\alpha[0,\ell_k]$ is a Banach algebra, each product
\[
        a u_{k,0},\qquad a u_{k,-},\qquad a u_{k,+}
\]
belongs to $H^\alpha[0,\ell_k]$, with the usual H\"older product estimate.  Substituting these products into \eqref{eq:Yrepresentation} gives a representation of $au$ in $Y_\alpha$ and yields the stated bound.
\end{proof}

\begin{remark}\label{rem:notalgebra}
The space $Y_\alpha$ is not an algebra in general.  The product of two functions containing logarithmic terms contains terms proportional to $(\log|t-t_j^d|)^2$.  This is not a defect for the present problem.  The equation requires multiplication of a logarithmic function by a piecewise H\"older coefficient, and this is precisely the module property of Lemma~\ref{lem:moduleproperty}.
\end{remark}

\begin{remark}[Vanishing leading coefficients versus strong H\"older regularity]\label{rem:vanishing_not_X}
The vanishing of all leading logarithmic coefficient functionals does not by itself imply membership in $X_\alpha$.  Indeed, on one arc with endpoint coordinate $s=0$, the function
\[
        u(s)=\chi(s)s^\alpha\log s
\]
belongs to $Y_\alpha$ and has zero leading logarithmic coefficient, because it may be represented with coefficient $u_-(s)=s^\alpha\in H^\alpha$ and $u_-(0)=0$.  Nevertheless,
\[
        \frac{|u(s)-u(0)|}{s^\alpha}=|\log s|\longrightarrow\infty,
\]
so $u\notin H^\alpha$ and hence $u\notin X_\alpha$.  Accordingly, throughout the paper the phrase \emph{log-free in the strong H\"older sense} means membership in $X_\alpha$, not merely the vanishing of the leading logarithmic coefficients.
\end{remark}

\section{The Cauchy singular integral equation in $X_\alpha\to Y_\alpha$}\label{sec:operator}

We study
\begin{equation}\label{eq:main}
        c(t)\varphi(t)+d(t)(S\varphi)(t)+(K\varphi)(t)=f(t),\qquad t\in\Ga\setminus\D,
\end{equation}
where
$
        (S\varphi)(t):=\frac{1}{\pi i}\PV\int_\Ga \frac{\varphi(\tau)}{\tau-t}\,d\tau
$
and
$
        (K\varphi)(t):=\frac{1}{2\pi i}\int_\Ga h(t,\tau)\varphi(\tau)\,d\tau.
$
The pointwise identity in \eqref{eq:main} is imposed only on $\Ga\setminus\D$.  Globally, the equation is interpreted as the equality
\[
        M\varphi=f\quad\hbox{in }Y_\alpha,
\]
which includes equality of the lateral logarithmic coefficients at every point of $\D$.  No finite point value of $S\varphi$ is required at a jump point when $\varphi$ has a nonzero jump.  The unknown is assumed to belong to $X_\alpha$, the coefficients satisfy
$
        c,d\in X_\alpha,
$
and the right-hand side belongs naturally to $Y_\alpha$.

\begin{assumption}[Regular kernel]\label{ass:kernel}
The kernel $h$ is continuous on $\Ga\times\Ga$ and is uniformly $\alpha$-H\"older in the first variable:
\[
        |h(t,\tau)-h(s,\tau)|\le C_h|t-s|^\alpha,
        \qquad t,s,\tau\in\Ga.
\]
\end{assumption}

\begin{lemma}[Mapping property of the regular part]\label{lem:Kmap}
Under Assumption~\ref{ass:kernel},
\[
        K:X_\beta\to H^\alpha(\Ga)\hookrightarrow X_\beta\hookrightarrow Y_\beta
\]
for every $0<\beta<\alpha<1$.  In particular, $K:X_\beta\to X_\beta$ is compact.
\end{lemma}

\begin{proof}
The boundedness in the sup norm follows from continuity of $h$ and boundedness of $\varphi$ on the finitely many arcs.  For $t,s\in\Ga$,
\[
        |K\varphi(t)-K\varphi(s)|
        \le C\int_\Ga |h(t,\tau)-h(s,\tau)| |\varphi(\tau)|\,|d\tau|
        \le C|t-s|^\alpha\|\varphi\|_{X_\beta}.
\]
Hence $K\varphi\in H^\alpha(\Ga)$.  The compactness into $X_\beta$ follows from the compact embedding $H^\alpha(\Ga)\hookrightarrow H^\beta(\Ga)$.
\end{proof}

The boundedness theory of singular integral operators in weighted and piecewise H\"older settings is closely related to the classical results of Duduchava~\cite{DuduchavaI,DuduchavaII}.

\begin{lemma}[Local logarithmic part of the Cauchy transform of a step]\label{lem:local_step_log}
Let $0<\alpha<1$, let $t_0\in\Ga$, and let $H$ be a left-continuous step function which is constant on the two sides of $t_0$ and has jump $J=[H]_{t_0}$.  Then, on each side of $t_0$,
\[
        (SH)(t)=\kappa^\pm J\log|t-t_0|+r^\pm(t),
        \qquad t\to t_0\pm0,
\]
where $\kappa^\pm\ne0$ depend only on the orientation of $\Ga$ and on the normalization of $S$, while $r^\pm\in H^\alpha$ up to the endpoint on the corresponding side.
\end{lemma}

\begin{proof}
It is enough to consider a sufficiently small neighbourhood of $t_0$.  Use a local parameter $s\in(-\varepsilon,\varepsilon)$ with $\gamma(0)=t_0$, $\gamma\in C^2$, and $\gamma'(0)\ne0$.  After subtracting the constant value on one side of the step, the only singular contribution is, up to the choice of the relevant half-arc,
\[
        I(x):=\frac{J}{\pi i}\int_0^\varepsilon
        \frac{\gamma'(s)}{\gamma(s)-\gamma(x)}\,ds,
        \qquad t=\gamma(x).
\]
On each side of $x=0$, choose logarithmic branches continuously on the subarcs that do not contain the pole.  If $x\in(0,\varepsilon)$, the integral is understood in the principal-value sense and the branches are chosen separately on $(0,x)$ and $(x,\varepsilon)$.  Since
\[
        \frac{d}{ds}\Log(\gamma(s)-\gamma(x))
        =\frac{\gamma'(s)}{\gamma(s)-\gamma(x)},
\]
the fundamental theorem of calculus, with the principal-value limit when necessary, yields
\[
        I(x)=-\frac{J}{\pi i}\Log(\gamma(0)-\gamma(x))+F^\pm(x),
\]
where $F^\pm$ is $C^1$ near $x=0$.  Indeed, apart from a side-dependent constant produced by the choice of logarithmic branches, $F^\pm$ is a constant multiple of $\Log(\gamma(\varepsilon)-\gamma(x))$, whose argument stays away from zero for sufficiently small $|x|$.

Now factor
\[
        \gamma(0)-\gamma(x)=-x\,q(x),
        \qquad
        q(x):=\int_0^1\gamma'(\theta x)\,d\theta.
\]
Because $\gamma\in C^2$, the function $q$ belongs to $C^1$, satisfies $q(0)=\gamma'(0)\ne0$, and is nonzero for sufficiently small $|x|$.  Hence, on each side of zero,
\[
        \Log(\gamma(0)-\gamma(x))
        =\log|x|+C^\pm+\Log q(x),
\]
where $\Log q\in C^1$.  Every $C^1$ function on a compact interval belongs to $H^\alpha$ for each $0<\alpha<1$.  Local equivalence of $|x|$ and $|t-t_0|$ therefore gives
\[
        (SH)(t)=\kappa^\pm J\log|t-t_0|+r^\pm(t),
\]
with $\kappa^\pm\ne0$ and $r^\pm\in H^\alpha$ on the corresponding closed side.  The part of the contour outside the local neighbourhood has a nonsingular kernel and contributes a $C^1$, hence $H^\alpha$, function.  This proves the assertion.
\end{proof}

\begin{proposition}[Cauchy operator from $X_\alpha$ to $Y_\alpha$]\label{prop:Smap}
For every $0<\alpha<1$, the Cauchy singular integral operator maps
\[
        S:X_\alpha\to Y_\alpha
\]
boundedly.
\end{proposition}

\begin{proof}
Let $\varphi\in X_\alpha$.  The relative Heaviside decomposition proved later in Proposition~\ref{prop:decomposition} gives
\[
        \varphi=\varphi_C+\sum_{j=2}^{n_d}\gamma_jG_j,
        \qquad \varphi_C\in H^\alpha(\Ga),
\]
with
$
        \|\varphi_C\|_{H^\alpha(\Ga)}+
        \sum_{j=2}^{n_d}|\gamma_j|
        \le C\|\varphi\|_{X_\alpha}.
$
The classical boundedness of the Cauchy singular operator on H\"older functions on a $C^2$ contour gives
$
        S\varphi_C\in H^\alpha(\Ga)\subset Y_\alpha.
$

It remains to consider the finitely many step functions $G_j$.  Each $G_j$ is constant on every arc of $\Ga\setminus\D$ and has only two jumps.  Applying Lemma~\ref{lem:local_step_log} at each jump point of $G_j$, and observing that the parts of the contour away from the jumps give H\"older contributions, we obtain on each side of every point $t_k^d$
\[
        (SG_j)(t)=A_{j,k}^{\pm}\log|t-t_k^d|+r_{j,k}^{\pm}(t),
\]
where $A_{j,k}^{\pm}$ are constants and the remainders $r_{j,k}^{\pm}$ are H\"older on the corresponding closed side.  Thus $SG_j\in Y_\alpha$.  Since there are finitely many such functions, the estimate follows by linearity.
\end{proof}

\begin{theorem}[Boundedness of the full operator]\label{thm:Mbounded}
Let $0<\alpha<1$, let $c,d\in X_\alpha$, and let Assumption~\ref{ass:kernel} hold.  Then
\[
        M:=cI+dS+K:X_\alpha\to Y_\alpha
\]
is bounded.
\end{theorem}

\begin{proof}
For $\varphi\in X_\alpha$, Lemma~\ref{lem:Xalgebra} gives $c\varphi\in X_\alpha\subset Y_\alpha$.  Proposition~\ref{prop:Smap} gives $S\varphi\in Y_\alpha$.  Since $d\in X_\alpha$, Lemma~\ref{lem:moduleproperty} gives $dS\varphi\in Y_\alpha$.  Finally, Lemma~\ref{lem:Kmap} gives $K\varphi\in Y_\alpha$.  Combining the estimates proves the result.
\end{proof}

\subsection{Well-posedness in the pair of spaces}\label{subsec:wellposedXY}

The mapping result in Theorem~\ref{thm:Mbounded} should be distinguished from invertibility.  In the discontinuous case considered here, the natural operator is
\[
        M:X_\beta\longrightarrow Y_\beta,
\]
not an operator from $X_\beta$ into itself.  Therefore the standard single-space Fredholm formulation based on the projectors $(I+S)/2$ and $(I-S)/2$ cannot be used without modification, because $S$ does not generally map $X_\beta$ into $X_\beta$ when the density has nonzero jumps.

Moreover, one should not expect $M(X_\beta)$ to coincide with all of $Y_\beta$.  The logarithmic coefficients of an element $M u$ are not arbitrary: they are generated by the jumps of $u$, by the coefficient $d$, and by the local logarithmic constants of the Cauchy operator.  Thus the appropriate continuous well-posedness condition for the present analysis is stability on the solution class, or equivalently bounded invertibility onto the range of $M$.

We shall use the following assumptions when the convergence theorem is proved:
\begin{enumerate}[label=(B\arabic*)]
\item $c,d\in X_\alpha$;
\item the regular kernel $h$ satisfies Assumption~\ref{ass:kernel};
\item the operator $M=cI+dS+K$ is bounded from $X_\beta$ to $Y_\beta$;
\item there exists a constant $C_M>0$ such that
\[
        \|u\|_{X_\beta}\le C_M\|Mu\|_{Y_\beta},
        \qquad u\in X_\beta;
\]
\item the right-hand side satisfies $f\in Y_\alpha\cap M(X_\alpha)$, and the exact solution belongs to $X_\alpha$.
\end{enumerate}
The fourth condition gives uniqueness and continuous dependence.  The fifth condition is a range condition: it says that the prescribed logarithmic right-hand side is compatible with the jump structure of a piecewise H\"older solution.  These assumptions replace the classical invertibility conditions in a single piecewise H\"older space.

\subsection{Compatibility when the right-hand side is log-free}\label{subsec:logfree}

The main analysis only requires
\[
        dS\varphi\in Y_\alpha,
\]
which follows from $S\varphi\in Y_\alpha$, $d\in X_\alpha$, and the module property $X_\alpha Y_\alpha\subset Y_\alpha$.  Thus the condition $dS\varphi\in X_\alpha$ is not an assumption of the collocation or convergence theory.  It becomes relevant only in the special case in which the prescribed right-hand side is required to belong to the smaller space $X_\alpha$.

Let $\varphi\in X_\alpha$ and suppose that $\varphi$ has a jump at $t_j^d$.  On the two sides of this point,
\begin{equation}\label{eq:localS}
        (S\varphi)(t)=A_j^\pm[\varphi]_{t_j^d}\log|t-t_j^d|+r_j^\pm(t),
\end{equation}
where $r_j^\pm$ is H\"older after subtraction of the displayed logarithmic term and $A_j^\pm\ne0$ are the local Cauchy constants determined by the orientation and normalization.  Since $c\varphi$ and $K\varphi$ belong to $X_\alpha$, every logarithmic term in $M\varphi$ is generated by $dS\varphi$.  Hence the leading logarithmic coefficient on either side of $t_j^d$ is proportional to
\[
        d(t_j^d\pm0)[\varphi]_{t_j^d}.
\]

\begin{proposition}[Necessary compatibility for $f\in X_\alpha$]\label{prop:compatibility}
Assume that $\varphi\in X_\alpha$ solves \eqref{eq:main} with $f\in X_\alpha$.  If $d(t_j^d+0)\ne0$ or $d(t_j^d-0)\ne0$, then
\[
        [\varphi]_{t_j^d}=0.
\]
Consequently, a nonzero jump of the solution at a point where $d$ does not vanish forces the right-hand side to contain a logarithmic term and therefore to lie naturally in $Y_\alpha\setminus X_\alpha$.
\end{proposition}

\begin{proof}
The lateral logarithmic coefficients of $f$, $c\varphi$, and $K\varphi$ vanish.  Formula \eqref{eq:localS} shows that the corresponding coefficients of $dS\varphi$ are nonzero constant multiples of $d(t_j^d\pm0)[\varphi]_{t_j^d}$.  Their vanishing therefore implies the stated conclusion.
\end{proof}

Pointwise vanishing of $d$ at the endpoint is not by itself sufficient for $dS\varphi\in X_\alpha$.  For example, in a one-sided coordinate $s\downarrow0$,
\[
        d(s)=s^\alpha\in H^\alpha,
        \qquad d(0)=0,
        \qquad s^\alpha\log s\notin H^\alpha.
\]
A sufficient strong cancellation condition on each side of $t_j^d$ is
\begin{equation}\label{eq:strong_log_cancellation}
        d(t_j^d\pm0)[\varphi]_{t_j^d}=0,
        \qquad
        [\varphi]_{t_j^d}\bigl(d(t)-d(t_j^d\pm0)\bigr)
        \log|t-t_j^d|\in H^\alpha.
\end{equation}
It holds, for example, if the jump vanishes, if $d$ is identically zero in a one-sided neighbourhood of the jump, or if $d(t_j^d\pm0)=0$ and $d$ has local regularity $H^\delta$ for some $\delta>\alpha$.  Otherwise the logarithmic target space $Y_\alpha$ and the separate logarithmic-coefficient equations used below are essential.

\section{Computational realization of the regularized B-spline--Heaviside collocation method}\label{sec:relative}\label{sec:computational}

This section converts the functional framework into an implementable scheme.  We first recall the relative Heaviside representation on a closed contour and the periodic B-spline approximation space, because these two ingredients determine the unknowns of the numerical method.  We then give the matrix-level realization of the regularized collocation equations.  The key point is that the logarithmic part generated by the Cauchy action on discontinuous trial functions is not sampled by ordinary point values; it is isolated through logarithmic-coefficient equations, while the remaining point equations are imposed only at nodes separated from the jump set.

\subsection{Relative Heaviside functions on a closed contour}

A step function on a closed contour cannot have a single nonzero jump.  Indeed, after one complete circuit around $\Ga$, the value of a single-valued function must return to its initial value; therefore the sum of all oriented jumps must be zero.  This observation motivates the following relative construction.

Fix the base jump point $t_1^d$.  For $j=2,\ldots,n_d$, let $\operatorname{arc}(a,b]$ denote the positively oriented arc that starts immediately after $a$ and ends at $b$.  We normalize the relative Heaviside function $G_j$ by
\begin{equation}\label{eq:relative_heaviside_definition}
        G_j(t):=
        \begin{cases}
        0, & t\in\operatorname{arc}(t_1^d,t_j^d],\\
        1, & t\in\operatorname{arc}(t_j^d,t_1^d].
        \end{cases}
\end{equation}
Thus $G_j$ is left-continuous at both endpoints, $G_j(t_j^d)=0$, and $G_j(t_1^d)=1$.  Its jumps are
\begin{equation*}
        [G_j]_{t_j^d}=1,
        \qquad [G_j]_{t_1^d}=-1,
        \qquad [G_j]_{t_i^d}=0\quad (i\ne1,j).
\end{equation*}
The normalization \eqref{eq:relative_heaviside_definition} removes the otherwise possible additive constant.  The functions $G_2,\ldots,G_{n_d}$ form a basis of the normalized step space consisting of step functions with jumps only at $\D$ that vanish on the first open arc leaving $t_1^d$.  The missing additive constant is absorbed into the continuous component; together with the continuous spline space, these functions give the direct sum used below.

\begin{proposition}[Relative decomposition]\label{prop:decomposition}
Let $0<\alpha<1$ and let $\varphi\in X_\alpha$.  Set
\[
        \gamma_j:=[\varphi]_{t_j^d},\qquad j=2,\ldots,n_d.
\]
Then
\begin{equation*}
        \varphi=\varphi_C+\sum_{j=2}^{n_d}\gamma_jG_j,
\end{equation*}
where $\varphi_C\in H^\alpha(\Ga)$.  Moreover,
\begin{equation}\label{eq:decompestimate}
        \|\varphi_C\|_{H^\alpha(\Ga)}+
        \sum_{j=2}^{n_d}|\gamma_j|
        \le C\|\varphi\|_{X_\alpha}.
\end{equation}
The representation is unique.
\end{proposition}

\begin{proof}
Since $\varphi$ is single-valued on the closed contour, the sum of all oriented jumps is zero:
\[
        \sum_{j=1}^{n_d}[\varphi]_{t_j^d}=0.
\]
Thus
$
        [\varphi]_{t_1^d}=-\sum_{j=2}^{n_d}\gamma_j.
$
Define
$
        \varphi_H:=\sum_{j=2}^{n_d}\gamma_jG_j,
        \qquad \varphi_C:=\varphi-\varphi_H.
$
At $t_j^d$, $j\ge2$, the jump of $\varphi_H$ is $\gamma_j$, and at $t_1^d$ it is $-\sum_{j=2}^{n_d}\gamma_j=[\varphi]_{t_1^d}$.  Hence $\varphi_C$ has zero jump at all points of $\D$ and extends continuously to the entire contour.

On each arc of $\Ga\setminus\D$, the function $\varphi_H$ is constant.  Therefore the local H\"older seminorm of $\varphi_C$ on that arc is the same as the local H\"older seminorm of $\varphi$.  To obtain global H\"older continuity, connect two arbitrary points by a finite chain of subarcs whose endpoints are points of $\D$.  Summing the local estimates and using equivalence of the arclength and Euclidean metrics gives a global $H^\alpha$ estimate.  The jump coefficients are bounded by $2\|\varphi\|_{\infty}$.  This proves \eqref{eq:decompestimate}.  Uniqueness follows because a continuous function cannot cancel a nonzero jump combination of the $G_j$.
\end{proof}

\subsection{Periodic B-spline spaces}

Let $m\ge2$ be fixed.  For $n_B\in\N$, set $h_B=2\pi/n_B$.  We choose a phase shift $\rho_{n_B}\in[0,h_B)$ and define parameter knots
\[
        \theta_k^{(n_B)}=\rho_{n_B}+(k-1)h_B,
        \qquad k=1,\ldots,n_B,
\]
with values understood modulo $2\pi$.  The corresponding contour knots are
\[
        t_k^B=\gamma(\theta_k^{(n_B)}),
        \qquad k=1,\ldots,n_B.
\]
Let $S_{n_B}$ denote the space of periodic splines of order $m$ on the parameter mesh, pulled back to $\Ga$ by $\gamma$.  Equivalently, if $\widetilde B_{m,k}$ are the standard periodic B-splines on the parameter circle, then
\[
        B_{m,k}(\gamma(\theta))=\widetilde B_{m,k}(\theta),
        \qquad k=1,\ldots,n_B,
\]
form a basis of $S_{n_B}$.

For the regularized collocation scheme used below, no separate mesh-separation theory is needed.  We only require that the spline nodes do not coincide with the jump points and that they do not approach them faster than the mesh size.  Thus, let $t_j^d=\gamma(\theta_j^d)$.  We impose the separation condition: there exists $\eta>0$, independent of $n_B$, such that
\begin{equation*}
        \dist(\theta_j^d-\rho_{n_B},h_B\mathbb Z)\ge \eta h_B,
        \qquad j=1,\ldots,n_d.
\end{equation*}
For every $n_B$ such a phase can be selected when $\eta<1/(2n_d)$, because the union of forbidden phase intervals has total length strictly smaller than $h_B$.  This condition ensures that all point evaluations used in the collocation equations are taken away from the logarithmic singularities.

We define the enriched B-spline--Heaviside trial space by
\begin{equation*}
        V_{n_B,H}:=S_{n_B}\oplus\spanop\{G_2,\ldots,G_{n_d}\}.
\end{equation*}
Thus
\[
        \dim V_{n_B,H}=n_B+n_d-1.
\]
Every $v\in V_{n_B,H}$ belongs to $X_\alpha$ for every $0<\alpha<1$.

\begin{proposition}[Approximation by enriched splines]\label{prop:approximation}
Let $0<\beta<\alpha<1$ and let $\varphi\in X_\alpha$.  Then there exists $v_{n_B}\in V_{n_B,H}$ such that
\begin{equation}\label{eq:approximation}
        \|\varphi-v_{n_B}\|_{X_\beta}
        \le C h_B^{\alpha-\beta}\|\varphi\|_{X_\alpha}.
\end{equation}
\end{proposition}

\begin{proof}
Use the decomposition
$
        \varphi=\varphi_C+\sum_{j=2}^{n_d}\gamma_jG_j,
         \varphi_C\in H^\alpha(\Ga).
$
We use the following standard Jackson-type estimate for periodic splines in the H\"older scale: for every $g\in H^\alpha(\Ga)$ there exists $s_{n_B}\in S_{n_B}$ such that
\begin{equation}\label{eq:periodic_spline_holder_estimate}
        \|g-s_{n_B}\|_{H^\beta(\Ga)}
        \le C h_B^{\alpha-\beta}\|g\|_{H^\alpha(\Ga)},
        \qquad 0<\beta<\alpha<1.
\end{equation}
For completeness, we recall why this estimate holds.  Pulling $g$ back by the smooth parametrization of $\Ga$, the problem becomes the corresponding periodic spline approximation problem on $[0,2\pi]$.  Let $Q_{n_B}$ be a standard local periodic B-spline quasi-interpolation operator of order $m$; see, for example, the classical spline approximation estimates in de Boor~\cite[Chs.~6--7]{deBoor} or Schumaker~\cite[Chs.~5--6]{Schumaker}.  Such operators are uniformly bounded in $H^\alpha$ and satisfy
\[
        \|g-Q_{n_B}g\|_{C(\Ga)}
        \le C h_B^\alpha [g]_{\alpha,\Ga},
        \qquad
        \|g-Q_{n_B}g\|_{H^\alpha(\Ga)}
        \le C\|g\|_{H^\alpha(\Ga)}.
\]
The elementary interpolation inequality
\[
        [w]_{\beta,\Ga}
        \le C\|w\|_{C(\Ga)}^{1-\beta/\alpha}
              [w]_{\alpha,\Ga}^{\beta/\alpha},
        \qquad 0<\beta<\alpha<1,
\]
applied to $w=g-Q_{n_B}g$ gives \eqref{eq:periodic_spline_holder_estimate}.  Applying \eqref{eq:periodic_spline_holder_estimate} to $g=\varphi_C$, we choose $s_{n_B}\in S_{n_B}$ such that
\[
        \|\varphi_C-s_{n_B}\|_{H^\beta(\Ga)}
        \le C h_B^{\alpha-\beta}\|\varphi_C\|_{H^\alpha(\Ga)}.
\]
Set
$
        v_{n_B}:=s_{n_B}+\sum_{j=2}^{n_d}\gamma_jG_j.
$
The jump part is represented exactly, and therefore
$
        \varphi-v_{n_B}=\varphi_C-s_{n_B}.
$
The $X_\beta$ norm of a globally $H^\beta$ function is bounded by its $H^\beta(\Ga)$ norm.  Using Proposition~\ref{prop:decomposition} gives \eqref{eq:approximation}.
\end{proof}

\begin{remark} \label{rem:weaker_norm}
Proposition~\ref{prop:approximation} shows that, in the presence of jump discontinuities, spline-based approximation is naturally formulated in a weaker piecewise H\"older norm.  For every $\varphi\in X_\alpha=PH^\alpha(\Ga,\D)$ and every $0<\beta<\alpha<1$, the enriched B-spline--Heaviside spaces produce approximants in $V_{n_B,H}$ that converge to $\varphi$ in the $X_\beta$-norm, with a rate determined by the H\"older regularity of the continuous component of $\varphi$.

This weaker norm is not merely a technical convenience.  The full space $X_\alpha$, equipped with the strong $\alpha$-H\"older norm on each continuity arc, is nonseparable.  Consequently, one should not expect a sequence of finite-dimensional spline-type spaces, or any countable approximation family, to be dense in $X_\alpha$ with respect to the $X_\alpha$-norm.  The use of $X_\beta$, $0<\beta<\alpha$, provides the appropriate approximation scale: it is strong enough to control the solution on each continuity arc, while it avoids the intrinsic nonseparability obstruction at the endpoint exponent.
\end{remark}

\subsection{Unknowns and block structure of the computed system}\label{subsec:blockstructure}

The numerical unknown is written in the form
\begin{equation}\label{eq:comp_ansatz}
        \varphi_{n_B}^H(t)=s_{n_B}(t)+g_H(t),
        \qquad
        s_{n_B}(t)=\sum_{k=1}^{n_B} a_k B_{m,k}(t),
        \qquad
        g_H(t)=\sum_{j=2}^{n_d}\gamma_jG_j(t).
\end{equation}
Here $a=(a_1,\ldots,a_{n_B})^T$ are the spline coefficients and
$\gamma=(\gamma_2,\ldots,\gamma_{n_d})^T$ are the independent jump amplitudes.  The relative Heaviside basis is essential on a closed contour: it automatically enforces the zero-total-jump constraint and avoids an artificial discontinuity at the periodic seam.

The action of the operator is decomposed as
\begin{equation*}
        M\varphi_{n_B}^H=M s_{n_B}+M g_H.
\end{equation*}
The spline part $s_{n_B}$ is continuous; hence $Ss_{n_B}$ is log-free.  Therefore all endpoint logarithmic terms in $M\varphi_{n_B}^H$ are generated by the Cauchy transforms of the relative Heaviside functions.  This observation gives a natural computational splitting:
\begin{enumerate}[label=\textup{(\alph*)}]
\item determine the Heaviside amplitudes from logarithmic-coefficient equations;
\item subtract their explicit logarithmic contribution from the point equations;
\item solve the resulting point-collocation system for the spline coefficients.
\end{enumerate}
The formal functional version of this construction is stated in Section~\ref{sec:collocation}.  The present subsection describes the actual matrix equations used in computation.

Let $\Lambda_2,\ldots,\Lambda_{n_d}$ be the selected logarithmic-coefficient functionals.  In the common implementation one chooses one lateral logarithmic coefficient at every non-base jump point, on a side where the corresponding lateral value of $d$ is nonzero.  The logarithmic equations are
\begin{equation*}
        \Lambda_i\!\left(M\varphi_{n_B}^H-f\right)=0,
        \qquad i=2,\ldots,n_d.
\end{equation*}
Since $c\varphi_{n_B}^H$, $K\varphi_{n_B}^H$, and $Ss_{n_B}$ are log-free, these equations reduce to the small system
\begin{equation}\label{eq:comp_Rgamma}
        R\gamma=b^{\log},
        \qquad
        R_{ij}=\Lambda_i(dSG_j),
        \qquad
        b_i^{\log}=\Lambda_i(f),
        \qquad i,j=2,\ldots,n_d.
\end{equation}
Thus the jump amplitudes are recovered before the main point-collocation system is assembled.  If $f$ is log-free, then $b^{\log}=0$; the equations impose cancellation of the logarithmic residual.  If $f\in Y_\alpha\setminus X_\alpha$, the vector $b^{\log}$ prescribes the logarithmic singular part of the right-hand side.

Let $\xi_i=\gamma(\vartheta_i)$, $i=1,\ldots,n_B$, be the point-collocation nodes, chosen away from $D$.  After $\gamma$ is known, the point equations
\begin{equation*}
        M\varphi_{n_B}^H(\xi_i)=f(\xi_i),
        \qquad i=1,\ldots,n_B,
\end{equation*}
become the dense linear system
\begin{equation}\label{eq:comp_Aa}
        A a=r,
\end{equation}
where
\begin{equation}\label{eq:comp_A_entries}
        A_{ik}=c(\xi_i)B_{m,k}(\xi_i)
        +d(\xi_i)(SB_{m,k})(\xi_i)
        +(KB_{m,k})(\xi_i)
\end{equation}
and
\begin{equation}\label{eq:comp_rhs_entries}
        r_i=f(\xi_i)-\sum_{j=2}^{n_d}\gamma_j
        \left[ c(\xi_i)G_j(\xi_i)
        +d(\xi_i)(SG_j)(\xi_i)
        +(KG_j)(\xi_i)\right].
\end{equation}
Equations \eqref{eq:comp_Rgamma} and \eqref{eq:comp_Aa} are the finite-dimensional computational counterpart of the regularized collocation equations.  In block form the unknown vector $(a,\gamma)$ satisfies an equivalent triangular system
\begin{equation*}
        \begin{pmatrix}
        A & H\\
        0 & R
        \end{pmatrix}
        \begin{pmatrix}a\\ \gamma\end{pmatrix}
        =
        \begin{pmatrix}f^p\\ b^{\log}\end{pmatrix},
\end{equation*}
where $f_i^p=f(\xi_i)$ and
\[
        H_{ij}=c(\xi_i)G_j(\xi_i)+d(\xi_i)(SG_j)(\xi_i)+(KG_j)(\xi_i).
\]
The triangular form is often preferable numerically because the small logarithmic block $R$ is assembled and solved first.

\subsection{Cauchy action on the spline part}\label{subsec:Sspline_comp}

For a continuous spline basis function $B_{m,k}$, the Cauchy singular integral is evaluated by singularity subtraction.  Using the normalization of $S$ and the identity
\begin{equation*}
        \frac{1}{\pi i}\,\PV\!\int_\Ga \frac{d\tau}{\tau-t}=1,
        \qquad t\in\Ga,
\end{equation*}
we write
\begin{equation}\label{eq:comp_S_spline}
        (SB_{m,k})(t)=B_{m,k}(t)
        +\frac{1}{\pi i}\int_\Ga
        \frac{B_{m,k}(\tau)-B_{m,k}(t)}{\tau-t}\,d\tau.
\end{equation}
The integrand in the second term has a removable singularity at $\tau=t$, because the numerator vanishes there.  In the parameter variable $t=\gamma(\theta)$ this becomes
\begin{equation*}
        (SB_{m,k})(\gamma(\theta_i))=B_{m,k}(\gamma(\theta_i))
        +\frac{1}{\pi i}\int_0^{2\pi}
        \frac{B_{m,k}(\gamma(\eta))-B_{m,k}(\gamma(\theta_i))}
        {\gamma(\eta)-\gamma(\theta_i)}\gamma'(\eta)\,d\eta.
\end{equation*}
At $\eta=\theta_i$ the integrand is replaced by its finite limiting value.  In the implementation the integral is split at spline knots, at jump parameters, and at the collocation parameter when needed.  On every resulting panel the integrand is regular, and a composite trapezoidal rule is applied.

\subsection{Cauchy action on the relative Heaviside part}\label{subsec:SHeaviside_comp}

The Cauchy transform of a relative Heaviside function is evaluated analytically, with the principal value treated explicitly.  For the normalization \eqref{eq:relative_heaviside_definition}, set
\[
        \Gamma_j^1:=\operatorname{arc}(t_j^d,t_1^d],
        \qquad
        \Gamma_j^0:=\operatorname{arc}(t_1^d,t_j^d].
\]
The function $G_j$ equals one on $\Gamma_j^1$ and zero on $\Gamma_j^0$.  At a collocation point $t\in\Ga\setminus\D$, exactly one of these two arcs contains $t$ in its interior.  Using
\[
        \frac{1}{\pi i}\PV\!\int_\Ga\frac{d\tau}{\tau-t}=1,
\]
we obtain the principal-value-safe formula
\begin{equation}\label{eq:comp_S_heaviside}
(SG_j)(t)=
\begin{cases}
\displaystyle
\frac{1}{\pi i}
\left[\Log_{\Gamma_j^1}(t_1^d-t)-\Log_{\Gamma_j^1}(t_j^d-t)\right],
& t\in\operatorname{int}\Gamma_j^0,\\[2.0ex]
\displaystyle
1-\frac{1}{\pi i}
\left[\Log_{\Gamma_j^0}(t_j^d-t)-\Log_{\Gamma_j^0}(t_1^d-t)\right],
& t\in\operatorname{int}\Gamma_j^1.
\end{cases}
\end{equation}
Here $\Log_{\Gamma_j^q}$ denotes a branch tracked continuously along the indicated oriented arc.  In each case the branch is followed only along an arc that does not contain the evaluation point, so the path $\tau-t$ never passes through zero.  In the second case the formula follows by subtracting the nonsingular integral over $\Gamma_j^0$ from the full-contour principal value.  Formula \eqref{eq:comp_S_heaviside} gives stable values of $(SG_j)(\xi_i)$ at all point-collocation nodes and also yields the endpoint logarithmic coefficients explicitly.  Thus the matrix $R$ in \eqref{eq:comp_Rgamma} is assembled without numerical differentiation and without evaluating limiting quotients near a singular endpoint.

\subsection{Regular integral part and quadrature}\label{subsec:K_comp}

The regular integral is evaluated in parameter form:
\begin{equation}\label{eq:comp_K_param}
        (Kv)(t)=\frac{1}{2\pi i}\int_0^{2\pi}
        h(t,\gamma(\eta))v(\gamma(\eta))\gamma'(\eta)\,d\eta.
\end{equation}
For $v=B_{m,k}$ the integrand is smooth on each spline element.  For $v=G_j$ the integral is split at the jump points so that the integrand is smooth on each continuity arc.  The same composite trapezoidal rule is then applied panel by panel.  In the numerical experiments the quadrature mesh is chosen finer than the spline mesh, so that the quadrature perturbation is smaller than the discretization error measured in the reported diagnostics.

\subsection{Implementation algorithm}\label{subsec:algorithm}

The complete computational procedure is summarized below.

\medskip
\noindent\textbf{Algorithm 1. Regularized B-spline--Heaviside collocation.}
\begin{enumerate}[label=\textbf{\arabic*.}, leftmargin=2.2em]
\item Choose the spline order $m$, the number of knots $n_B$, and a phase-shifted periodic mesh separated from the jump parameters.
\item Construct the periodic B-spline basis $B_{m,k}$, $k=1,\ldots,n_B$.
\item Construct the relative Heaviside functions $G_j$, $j=2,\ldots,n_d$.
\item Select logarithmic-coefficient functionals $\Lambda_2,\ldots,\Lambda_{n_d}$ and assemble
      $R_{ij}=\Lambda_i(dSG_j)$ and $b_i^{\log}=\Lambda_i(f)$.
\item Solve the small logarithmic system $R\gamma=b^{\log}$ for the Heaviside amplitudes.
\item For each collocation node $\xi_i$, compute $(SB_{m,k})(\xi_i)$ by \eqref{eq:comp_S_spline}, $(SG_j)(\xi_i)$ by \eqref{eq:comp_S_heaviside}, and the regular terms $(KB_{m,k})(\xi_i)$ and $(KG_j)(\xi_i)$ by \eqref{eq:comp_K_param}.
\item Assemble $A$ and $r$ from \eqref{eq:comp_A_entries}--\eqref{eq:comp_rhs_entries}.
\item Solve $Aa=r$.
\item Reconstruct $\varphi_{n_B}^H$ from \eqref{eq:comp_ansatz}; report the jump amplitudes, point residuals, logarithmic residuals, and error diagnostics when an exact solution is available.
\end{enumerate}
\medskip

The logarithmic system has size $n_d-1$, which is fixed in the typical situation where the number of prescribed jumps does not grow with $n_B$.  The point-collocation matrix is dense because of the Cauchy and Fredholm integral terms.  Direct solution by Gaussian elimination has cost $O(n_B^3)$, while assembly depends on the quadrature strategy and on whether the local support of the splines is exploited in evaluating the regular integral part.  For the sizes used in Section~\ref{sec:numerical}, dense direct solvers are sufficient.

\section{Regularized B-spline--Heaviside collocation}\label{sec:collocation}

The equation is interpreted in the target space $Y_\beta$.  This has an immediate consequence for collocation.  Ordinary point evaluation is a bounded operation on $X_\beta$, and it is also harmless for functions in $Y_\beta$ at points separated from the jump set.  However, a general element of $Y_\beta$ may have a logarithmic singularity at a point of $\D$ and therefore need not have a finite value there.  Thus one should not impose the equation at the jump points by the literal conditions
\[
        (M v)(t_j^d)=f(t_j^d).
\]
Instead, the endpoint information must be imposed through bounded functionals on $Y_\beta$.  In this section these functionals are chosen to be logarithmic-coefficient functionals.  They extract the coefficients of the logarithmic terms generated by the jumps of the trial function.  This is the part of the residual that is invisible to ordinary values away from the jump point but is essential for determining the Heaviside amplitudes.

\subsection{Logarithmic coefficients at a jump point}

Recall the notation introduced in Section~\ref{sec:spaces}.  If $u\in Y_\beta$ and if $t_j^d$ is a jump point, then on the two arcs adjacent to $t_j^d$ the function $u$ has local representations of the form
\begin{align*}
        u(t)&=a_j^+(t)\log |t-t_j^d|+r_j^+(t),
        &&t\to t_j^d+0,\\
        u(t)&=a_j^-(t)\log |t-t_j^d|+r_j^-(t),
        &&t\to t_j^d-0,
\end{align*}
where $a_j^\pm$ and $r_j^\pm$ are H\"older functions up to the corresponding endpoint.  The leading logarithmic coefficients are
\[
        \lambda_j^+(u):=a_j^+(t_j^d+0),\qquad
        \lambda_j^-(u):=a_j^-(t_j^d-0).
\]
They are independent of the particular cut-off representation of $u$ and define bounded linear functionals on $Y_\beta$.  Equivalently,
\[
        \lambda_j^\pm(u)
        =\lim_{t\to t_j^d\pm0}\frac{u(t)}{\log |t-t_j^d|},
\]
provided the remaining part is understood in the finite-part sense.  If $u\in X_\beta$, then $u$ is log-free and therefore
\[
        \lambda_j^+(u)=\lambda_j^-(u)=0,
        \qquad j=1,\ldots,n_d.
\]

For $u\in X_\beta$, the logarithmic part of $Su$ is determined by the jumps of $u$.  More precisely, for each side of each point $t_j^d$ there are nonzero constants $\kappa_j^+$ and $\kappa_j^-$, depending only on the orientation of $\Gamma$ and on the normalization of $S$, such that
\begin{equation}\label{eq:local_log_coeff}
        \lambda_j^\pm(Su)=\kappa_j^\pm [u]_{t_j^d},
        \qquad u\in X_\beta.
\end{equation}
The exact numerical value of $\kappa_j^\pm$ is immaterial for the abstract analysis.  What matters is that these constants are nonzero.  Relation~\eqref{eq:local_log_coeff} follows from the local model of the Cauchy principal value integral of a unit step.  The difference between the contour kernel and the model kernel produces only a H\"older remainder, hence no additional logarithmic coefficient.

Since $c u$ and $K u$ are log-free whenever $u\in X_\beta$, the logarithmic coefficients of $Mu=c u+dSu+Ku$ are produced only by $dSu$.  Thus
\begin{equation}\label{eq:M_log_coeff}
        \lambda_j^\pm(Mu)=d(t_j^d\pm0)\,\kappa_j^\pm [u]_{t_j^d},
        \qquad u\in X_\beta.
\end{equation}
This formula explains why logarithmic endpoint equations are the natural equations for the jump amplitudes.

\subsection{Choice of logarithmic observation functionals}

Let
\[
        \Xi_{n_B}=\{t_1^B,\ldots,t_{n_B}^B\}
\]
be the phase-shifted knot set.  By construction, $\Xi_{n_B}\cap\D=\emptyset$ and the knots are separated from $\D$.  For $u\in Y_\beta$, define point evaluation functionals
\[
        e_k(u):=u(t_k^B),\qquad k=1,\ldots,n_B.
\]
The values are finite because the mesh nodes do not coincide with the logarithmic singularities.

At the jump points we use $n_d-1$ logarithmic-coefficient functionals
\[
        \Lambda_2,\ldots,\Lambda_{n_d}\in Y_\beta^*.
\]
Each $\Lambda_j$ is a fixed linear combination of the lateral coefficient functionals $\lambda_i^\pm$.  In the simplest and most common case one takes
\begin{equation*}
        \Lambda_j=\lambda_j^{\sigma_j},
        \qquad \sigma_j\in\{+,-\},\quad j=2,\ldots,n_d,
\end{equation*}
choosing, at each non-base jump point $t_j^d$, one of the two lateral sides.  More generally, linear combinations may be useful if both sides are involved or if one wants to normalize the equations.

The purpose of the functionals $\Lambda_j$ is not to approximate point values at the jump points.  Their purpose is to observe the logarithmic part of the residual.  Since the relative Heaviside function $G_j$ has jump $+1$ at $t_j^d$ and jump $-1$ at the base point $t_1^d$, the action of $M$ on $G_j$ generates logarithmic terms at these two points.  The $n_d-1$ functionals $\Lambda_2,\ldots,\Lambda_{n_d}$ are chosen so that the logarithmic components generated by $G_2,\ldots,G_{n_d}$ can be distinguished.

\begin{assumption}[Logarithmic unisolvence]\label{ass:logunisolvence}
The matrix
\begin{equation*}
        \mathcal L=\bigl(\Lambda_i(MG_j)\bigr)_{i,j=2}^{n_d}
\end{equation*}
is nonsingular.
\end{assumption}

Assumption~\ref{ass:logunisolvence} means that the logarithmic endpoint equations determine the unknown jump amplitudes.  Indeed, if
\[
        v=s+\sum_{j=2}^{n_d}\gamma_jG_j\in V_{n_B,H},
\]
then the logarithmic part of $Mv$ is independent of the smooth spline component $s$ and is governed by the finite vector $\gamma=(\gamma_2,\ldots,\gamma_{n_d})$.  The matrix $\mathcal L$ is precisely the finite matrix that maps these jump amplitudes to the observed logarithmic coefficients.  Thus Assumption~\ref{ass:logunisolvence} is checked directly: after the functionals $\Lambda_j$ are chosen, one computes the finite matrix $\mathcal L$ of logarithmic coefficients and verifies that its determinant is nonzero.  No limiting process or infinite-dimensional argument is involved in this check.

\subsection{A verifiable sufficient condition for Assumption~\ref{ass:logunisolvence}}

Assumption~\ref{ass:logunisolvence} is abstract, but in the present setting it has a simple verifiable sufficient condition.  Suppose that, for each $j=2,\ldots,n_d$, one can choose a side $\sigma_j\in\{+,-\}$ such that
\begin{equation}\label{eq:d_nonzero_side}
        d(t_j^d\sigma_j0)\ne0.
\end{equation}
Here $d(t_j^d+0)$ and $d(t_j^d-0)$ denote the one-sided limits of $d$ at $t_j^d$.  Choose
\[
        \Lambda_j:=\lambda_j^{\sigma_j},
        \qquad j=2,\ldots,n_d.
\]
Then, using \eqref{eq:M_log_coeff} and the defining jumps of the relative Heaviside functions,
\begin{equation*}
        \Lambda_i(MG_j)
        =d(t_i^d\sigma_i0)\,\kappa_i^{\sigma_i}\,[G_j]_{t_i^d}
        =d(t_i^d\sigma_i0)\,\kappa_i^{\sigma_i}\,\delta_{ij},
        \qquad i,j=2,\ldots,n_d.
\end{equation*}
Since the constants $\kappa_i^{\sigma_i}$ are nonzero and the multipliers $d(t_i^d\sigma_i0)$ are nonzero by \eqref{eq:d_nonzero_side}, the matrix $\mathcal L$ is diagonal with nonzero diagonal entries.  Hence it is nonsingular.

Thus a practical sufficient condition is:
\begin{quote}
For every non-base jump point $t_j^d$, at least one lateral value of $d$ at $t_j^d$ is nonzero; then choose the logarithmic equation on such a side.
\end{quote}
If this condition fails at some jump point, then the Cauchy logarithm produced by a jump may be annihilated by the coefficient $d$ on both sides.  In that case the jump amplitude is not detected by logarithmic coefficients at that point, and one must either use a different set of endpoint functionals, such as finite-part functionals, or impose an additional compatibility equation.  The base point $t_1^d$ is not lost in this construction: because the relative Heaviside functions have zero total jump, the logarithmic coefficient at the base point is determined by the same amplitudes $\gamma_2,\ldots,\gamma_{n_d}$.  One may use the base-point logarithmic coefficient instead of one of the non-base coefficients if this gives a better-conditioned finite matrix, but one should still impose only $n_d-1$ independent endpoint equations.

For the relative Heaviside basis, the logarithmic block can be assembled directly.  If the endpoint functional is $\Lambda_i=\lambda_i^{\sigma_i}$, then
$
        \Lambda_i(dSG_j)
        =d(t_i^d\sigma_i0)\,\lambda_i^{\sigma_i}(SG_j).
$
With the unnormalized logarithmic coefficient, one has
$
        \lambda_i^{\sigma_i}(SG_j)
        =\kappa_i^{\sigma_i}[G_j]_{t_i^d}.
$
In the MATLAB implementation used for the experiments, the logarithmic coefficient functionals are normalized so that
$
        \lambda_i^{\sigma_i}(Sv)=[v]_{t_i^d},
         v\in X_\beta .
$
Therefore
\[
        \Lambda_i(dSG_j)=d(t_i^d\sigma_i0)[G_j]_{t_i^d}.
\]
Choosing the plus side at the non-base jump points $t_i^d$, $i=2,\ldots,n_d$, and using $[G_j]_{t_i^d}=\delta_{ij}$, gives
\[
        R_{ij}=d(t_i^d+0)\delta_{ij},
        \qquad i,j=2,\ldots,n_d.
\]
This is exactly the block assembled in the code.  It is nonsingular as long as the selected lateral values $d(t_i^d+0)$ are nonzero.  For the manufactured tests, $\Lambda_i(f)=d(t_i^d+0)\gamma_i^{\rm exact}$, which explains why the numerical jump error is at the level of machine precision.

\subsection{The observation operator and the scaled discrete residual norm}

Define
\begin{equation*}
        E_{n_B}:Y_\beta\to\C^{n_B+n_d-1}
\end{equation*}
by
\begin{equation*}
        E_{n_B}u:=
        \bigl(e_1(u),\ldots,e_{n_B}(u),
        \Lambda_2(u),\ldots,\Lambda_{n_d}(u)\bigr).
\end{equation*}
The unscaled algebraic residual norm is
\begin{equation*}
        \|E_{n_B}u\|_{Z_{n_B}}
        :=\max\left\{
        \max_{1\le k\le n_B}|u(t_k^B)|,
        \max_{2\le j\le n_d}|\Lambda_j(u)|
        \right\}.
\end{equation*}
Thus
\begin{equation*}
        E_{n_B}M v
        =\bigl((Mv)(t_1^B),\ldots,(Mv)(t_{n_B}^B),
        \Lambda_2(Mv),\ldots,\Lambda_{n_d}(Mv)\bigr).
\end{equation*}
The residual is measured by ordinary values away from $\D$ and by logarithmic coefficients at $\D$.  This is the regularized analogue of classical point collocation.

For the stability proof it is essential to use the eliminated point residual.  Let
\[
        v=s+\sum_{j=2}^{n_d}\gamma_jG_j\in V_{n_B,H},
        \qquad s\in S_{n_B},
\]
and set
\[
        \ell(v):=\bigl(\Lambda_2(Mv),\ldots,\Lambda_{n_d}(Mv)\bigr)^T .
\]
Since the spline component is continuous, it has no logarithmic coefficients; hence
\[
        \ell(v)=R\gamma,
        \qquad R=(\Lambda_i(MG_j))_{i,j=2}^{n_d} .
\]
Let
\[
        p_{n_B}(v):=\bigl((Mv)(t_1^B),\ldots,(Mv)(t_{n_B}^B)\bigr)^T
\]
and let
\[
        H_{n_B}:=\bigl((MG_j)(t_k^B)\bigr)_{\substack{k=1,\ldots,n_B\\ j=2,\ldots,n_d}}
\]
be the point-collocation coupling generated by the relative Heaviside functions.  If $R$ is nonsingular, define the eliminated, or regularized, point residual by
\begin{equation*}
        p_{n_B}^{\rm reg}(v)
        :=p_{n_B}(v)-H_{n_B}R^{-1}\ell(v).
\end{equation*}
For an exact trial function decomposition this is simply
\begin{equation*}
        p_{n_B}^{\rm reg}(v)
        =\bigl((Ms)(t_1^B),\ldots,(Ms)(t_{n_B}^B)\bigr)^T .
\end{equation*}
Thus the logarithmic equations determine the Heaviside amplitudes, and the eliminated point block measures only the log-free spline component.

For a compatible right-hand side $f\in Y_\beta$, the same elimination is used:
\[
        \ell(f):=\bigl(\Lambda_2(f),\ldots,\Lambda_{n_d}(f)\bigr)^T,
        \qquad
        p_{n_B}^{\rm reg}(f):=\bigl(f(t_1^B),\ldots,f(t_{n_B}^B)\bigr)^T
        -H_{n_B}R^{-1}\ell(f).
\]
We write
\[
        \widehat E_{n_B}Mv:=\bigl(p_{n_B}^{\rm reg}(v),\ell(v)\bigr),
        \qquad
        \widehat E_{n_B}f:=\bigl(p_{n_B}^{\rm reg}(f),\ell(f)\bigr).
\]
The original equations $E_{n_B}Mv=E_{n_B}f$ and the eliminated equations
\[
        \widehat E_{n_B}Mv=\widehat E_{n_B}f
\]
are equivalent whenever $R$ is nonsingular.

The stability estimate in the $X_\beta$-norm cannot be uniform with respect to the unscaled maximum norm of point residuals alone.  Even for the identity operator on a continuous spline space, the $H^\beta$ seminorm of a spline can grow like $h_B^{-\beta}$ times its nodal maximum.  The correct residual norm for the eliminated observation is therefore the scaled norm
\begin{equation}\label{eq:scaled_residual_norm}
        \|\widehat E_{n_B}Mv\|_{\widehat Z_{n_B}^{(\beta)}}
        :=\max\left\{
        h_B^{-\beta}\|p_{n_B}^{\rm reg}(v)\|_{\ell^\infty},
        \|\ell(v)\|_{\ell^\infty}
        \right\}.
\end{equation}
This scaling is the finite-dimensional counterpart of the inverse inequality for splines and is the form used in the stability and perturbation estimates below.

For general logarithmic $u\in Y_\beta$, the point values $u(t_k^B)$ near $\D$ may grow like $|\log h_B|$.  The convergence proof below does not require a uniform bound of $E_{n_B}$ on all of $Y_\beta$.  It uses the eliminated observation and applies point evaluations only to log-free residuals.

\begin{lemma}[Observation of log-free functions]\label{lem:observationX}
There exists $C>0$, independent of $n_B$, such that for all $u\in X_\beta$,
\[
        \|E_{n_B}u\|_{Z_{n_B}}
        \le C\|u\|_{X_\beta}.
\]
\end{lemma}

\begin{proof}
For $u\in X_\beta$, the point evaluations are bounded by the sup-norm part of $\|u\|_{X_\beta}$.  Since $u$ has no logarithmic terms, all logarithmic coefficient functionals vanish.  Therefore the second block of $E_{n_B}u$ is zero and the first block is uniformly controlled by $\|u\|_{X_\beta}$.
\end{proof}

\subsection{Exact regularized collocation scheme}

The exact regularized collocation problem is: find $\varphi_{n_B}^H\in V_{n_B,H}$ such that
\begin{equation}\label{eq:scheme}
        E_{n_B}M\varphi_{n_B}^H=E_{n_B}f.
\end{equation}
Equivalently, the equations consist of
\begin{align}
        (M\varphi_{n_B}^H)(t_k^B)&=f(t_k^B),
        &&k=1,\ldots,n_B,\notag\\
        \Lambda_j(M\varphi_{n_B}^H)&=\Lambda_j(f),
        &&j=2,\ldots,n_d.\label{eq:logcoll}
\end{align}
No value of $f$ at a logarithmic singularity is required.  If $f\in X_\alpha$, then $\Lambda_j(f)=0$, and the equations \eqref{eq:logcoll} impose the vanishing of the corresponding logarithmic coefficients of the residual.  If $f\in Y_\alpha\setminus X_\alpha$, the right-hand side of \eqref{eq:logcoll} prescribes exactly the logarithmic singular part that must be reproduced by the discontinuous component of the approximate solution.

\begin{remark}[Finite-part variants]\label{rem:finitepart}
Logarithmic coefficient equations are not the only possible regularized endpoint equations.  One may instead use finite-part endpoint equations after subtracting the logarithmic term explicitly.  Suppose, on a fixed side of $t_j^d$, that
\[
        u(t)=a(t)\log|t-t_j^d|+r(t),
\]
where $a$ and $r$ are H\"older up to the endpoint.  The logarithmic coefficient is $a(t_j^d\pm0)$, while the finite part is
\[
        \operatorname{f.p.}_{j}^{\pm}u
        :=\lim_{t\to t_j^d\pm0}
        \bigl(u(t)-a(t_j^d\pm0)\log|t-t_j^d|\bigr)
        =r(t_j^d\pm0).
\]
This functional is finite and bounded on $Y_\beta$ once the logarithmic coefficient is known.  A collocation method may therefore replace some logarithmic coefficient equations by finite-part equations, or use both types when an overdetermined least-squares formulation is desired.

The present paper uses logarithmic coefficient equations because they directly isolate the compatibility conditions created by jumps.  They are especially transparent for determining the Heaviside amplitudes: the logarithmic coefficient of $MG_j$ is given explicitly by \eqref{eq:M_log_coeff}.  Finite-part equations are useful in implementation when one wants to match the regular endpoint value after the singular term has been removed.  If finite-part functionals are used instead of the $\Lambda_j$, Assumption~\ref{ass:logunisolvence} must be replaced by the analogous nonsingularity of the finite matrix obtained from those endpoint functionals.
\end{remark}

\section{Conditional stability and convergence for exact operator action}\label{sec:exactconv}

We now prove a conditional convergence theorem for the exact collocation scheme.  The result is stated for a fixed pair $0<\beta<\alpha<1$.  The mesh-uniform discrete stability and scaled consistency hypotheses below are independent assumptions on the concrete discretization; they are not consequences of continuous stability alone.

\begin{assumption}[Continuous stability]\label{ass:continuous}
The operator $M:X_\beta\to Y_\beta$ is stable on its range: there exists $C_M>0$ such that
\begin{equation}\label{eq:continuous}
        \|u\|_{X_\beta}\le C_M\|Mu\|_{Y_\beta},
        \qquad u\in X_\beta.
\end{equation}
Equivalently, $M$ is injective and has a bounded inverse from the range $M(X_\beta)$, equipped with the $Y_\beta$-norm, onto $X_\beta$.
\end{assumption}

Estimate~\eqref{eq:continuous} is the continuous well-posedness estimate used in this paper.  It should not be read as the statement that $M$ maps $X_\beta$ onto the whole space $Y_\beta$.  In general, $M(X_\beta)$ is a proper subspace of $Y_\beta$, because the logarithmic coefficients of $Mu$ must be compatible with the jumps of $u$ and with the coefficient $d$.  Instead, \eqref{eq:continuous} says that the residual in the natural target norm controls the unknown in the solution norm.  A standard sufficient condition is that $M:X_\beta\to Y_\beta$ have closed range and trivial kernel; then the bounded inverse theorem gives \eqref{eq:continuous}.  In applications this is typically verified by separating the finitely many logarithmic coefficient equations from the regular equation for the continuous component.
A more concrete sufficient route to Assumption~6.1 is obtained after
separating the jump amplitudes from the continuous component. Using the
relative Heaviside decomposition
\[
u=u_C+\sum_{j=2}^{n_d}\gamma_jG_j,
\qquad u_C\in H^\beta(\Gamma),
\]
and applying the logarithmic-coefficient functionals, the equation can be
represented, at the continuous level, as a coupled block system
\[
\begin{pmatrix}
T & B\\
0 & R
\end{pmatrix}
\binom{u_C}{\gamma}
=
\binom{g}{\ell}.
\]
Here \(R\) is the finite logarithmic block and \(T\) is the log-free operator
acting on the continuous component after the logarithmic part has been
separated. If the relative Heaviside decomposition is used as a topological
identification of \(X_\beta\) with
\(H^\beta(\Gamma)\times\mathbb C^{n_d-1}\), if \(R\) is nonsingular, and if
\(T\) is Fredholm of index zero with trivial kernel, then \(T\) is boundedly
invertible and the block operator is boundedly invertible onto its range.
Consequently, the continuous stability estimate in Assumption~6.1 follows.
In the present paper we keep Assumption~6.1 as the abstract well-posedness
condition, because the Fredholm verification of \(T\) depends on the
coefficients and on the contour.

\begin{assumption}[Scaled uniform discrete stability]\label{ass:discrete}
There exist $C_s>0$ and $n_0\in\N$ such that, for all $n_B\ge n_0$ and all $v\in V_{n_B,H}$,
\begin{equation*}
        \|v\|_{X_\beta}
        \le C_s\|\widehat E_{n_B}Mv\|_{\widehat Z_{n_B}^{(\beta)}}.
\end{equation*}
\end{assumption}

Assumption~\ref{ass:discrete} is the discrete counterpart of continuous stability for the eliminated regularized observation.  The scaling in \eqref{eq:scaled_residual_norm} is unavoidable when the left-hand side is the $X_\beta$-norm: point values alone control the maximum norm of the spline component, while the $\beta$-H\"older seminorm is recovered by an inverse inequality.  The following proposition gives a concrete sufficient criterion which is directly connected with the matrices assembled in the implementation.

\begin{proposition}[Sufficient criterion for scaled discrete stability]\label{prop:block_scaled}
Let
\[
        v=s+g_\gamma,
        \qquad
        s\in S_{n_B},
        \qquad
        g_\gamma=\sum_{j=2}^{n_d}\gamma_jG_j.
\]
Assume that the following estimates hold uniformly for all sufficiently large $n_B$:
\begin{enumerate}[label=(S\arabic*)]
\item the logarithmic block is uniformly invertible,
\[
        \|\gamma\|_{\ell^\infty}
        \le C_R\|R\gamma\|_{\ell^\infty},
        \qquad \gamma\in\C^{n_d-1};
\]
\item the zero-jump spline point block is stable in the maximum norm,
\[
        \|s\|_{C(\Gamma)}
        \le C_A\max_{1\le k\le n_B}|(Ms)(t_k^B)|,
        \qquad s\in S_{n_B};
\]
\item the fixed relative Heaviside basis is uniformly bounded in $X_\beta$,
\[
        \|g_\gamma\|_{X_\beta}\le C_G\|\gamma\|_{\ell^\infty}.
\]
\end{enumerate}
Then Assumption~\ref{ass:discrete} holds.
\end{proposition}

\begin{proof}
Since $s$ is continuous, it has no jump contribution to the logarithmic block.  Therefore
\[
        \ell(v)=R\gamma,
\]
and (S1) gives
\[
        \|\gamma\|_{\ell^\infty}
        \le C_R\|\ell(v)\|_{\ell^\infty}.
\]
By the definition of the eliminated point residual,
\[
        p_{n_B}^{\rm reg}(v)=
        \bigl((Ms)(t_1^B),\ldots,(Ms)(t_{n_B}^B)\bigr)^T .
\]
Hence (S2) implies
\[
        \|s\|_{C(\Gamma)}
        \le C_A\|p_{n_B}^{\rm reg}(v)\|_{\ell^\infty}.
\]
For periodic splines of fixed order on a quasi-uniform mesh, the standard inverse inequality gives
$
        \|s\|_{H^\beta(\Gamma)}
        \le C_I h_B^{-\beta}\|s\|_{C(\Gamma)}.
$
Thus
\[
        \|s\|_{X_\beta}
        \le C h_B^{-\beta}\|p_{n_B}^{\rm reg}(v)\|_{\ell^\infty}.
\]
Using (S3) and the bound for $\gamma$, we obtain
\[
\begin{aligned}
        \|v\|_{X_\beta}
        &\le \|s\|_{X_\beta}+\|g_\gamma\|_{X_\beta}
        &\le C h_B^{-\beta}\|p_{n_B}^{\rm reg}(v)\|_{\ell^\infty}
             + C\|\ell(v)\|_{\ell^\infty}
        &\le C\|\widehat E_{n_B}Mv\|_{\widehat Z_{n_B}^{(\beta)}} .
\end{aligned}
\]
The constant is independent of $n_B$.
\end{proof}

\begin{remark}[How the criterion is checked]\label{rem:scaled_check}
Condition (S1) is completely explicit in the usual lateral normalization of Section~\ref{sec:collocation}: if $\lambda_i^{\sigma_i}(Sv)=[v]_{t_i^d}$ and $\Lambda_i=\lambda_i^{\sigma_i}$, then
\[
        R=\diag\bigl(d(t_2^d\sigma_2 0),\ldots,d(t_{n_d}^d\sigma_{n_d}0)\bigr).
\]
Thus (S1) follows from a positive lower bound for the selected lateral values of $d$.  Condition (S3) is automatic because the number of relative Heaviside functions is fixed.  The only genuinely discretization-dependent part is (S2), the usual stability of the zero-jump spline collocation block in the maximum norm.  For a fixed spline order, phase rule, and coefficient set, it is equivalent to a uniform lower bound for the inverse of the scaled spline collocation matrices.  The algebraic conditioning data in Section~\ref{subsec:num_conditioning} and the direct inverse-norm proxy in Section~\ref{subsec:num_scaled_diagnostics} provide finite-dimensional checks for the reported computations, but they do not establish a mesh-uniform bound as $n_B\to\infty$.
\end{remark}

\begin{assumption}[Scaled consistency of exact-jump approximants]\label{ass:consistency}
For the exact solution $\varphi\in X_\alpha$, there exist approximants $w_{n_B}\in V_{n_B,H}$ with the same jumps as $\varphi$ such that
\begin{equation}\label{eq:scaled_consistency}
        \|\varphi-w_{n_B}\|_{X_\beta}
        +\|\widehat E_{n_B}M(\varphi-w_{n_B})\|_{\widehat Z_{n_B}^{(\beta)}}
        \le C h_B^{\alpha-\beta}\|\varphi\|_{X_\alpha}.
\end{equation}
\end{assumption}

The first term in \eqref{eq:scaled_consistency} is the approximation property of Proposition~\ref{prop:approximation}.  The second term is an additional, strictly stronger collocation consistency requirement in the scaled residual norm.  It does not follow from the first term or from Lemma~\ref{lem:regularresidual}.  Indeed, because the point residual is multiplied by $h_B^{-\beta}$, the second term requires, in particular,
\[
        \max_{1\le k\le n_B}
        \bigl|M(\varphi-w_{n_B})(t_k^B)\bigr|
        =O(h_B^\alpha),
\]
whereas a mere $X_\beta$ bound of order $O(h_B^{\alpha-\beta})$ would not suffice after scaling.  For exact-jump approximants the logarithmic part of $M(\varphi-w_{n_B})$ vanishes, and the remaining estimate must be proved for the chosen spline approximant, coefficient set, and collocation nodes, or verified independently for the concrete scheme.

\begin{lemma}[Regular residual for exact-jump approximants]\label{lem:regularresidual}
Let $w_{n_B}\in V_{n_B,H}$ have the same jumps as the exact solution $\varphi\in X_\alpha$.  Then
\[
        r_{n_B}:=\varphi-w_{n_B}\in H^\beta(\Gamma)
\]
and
\begin{equation}\label{eq:regularresidual}
        Mr_{n_B}\in X_\beta,
        \qquad
        \|Mr_{n_B}\|_{X_\beta}\le C\|r_{n_B}\|_{X_\beta}.
\end{equation}
\end{lemma}

\begin{proof}
Since $w_{n_B}$ and $\varphi$ have the same jumps, $r_{n_B}$ has no jumps.  Hence $r_{n_B}\in H^\beta(\Gamma)$.  The Cauchy singular integral maps $H^\beta(\Gamma)$ boundedly into itself.  Since $c,d\in X_\alpha\subset X_\beta$, the products $cr_{n_B}$ and $dSr_{n_B}$ belong to $X_\beta$.  The regular part $Kr_{n_B}$ belongs to $H^\alpha(\Gamma)\subset X_\beta$.  This proves \eqref{eq:regularresidual}.
\end{proof}

\begin{theorem}[Conditional existence, uniqueness, and convergence]\label{thm:main}
Let $0<\beta<\alpha<1$.  Assume that $c,d\in X_\alpha$, that Assumption~\ref{ass:kernel} holds, and that $f\in Y_\alpha\cap M(X_\alpha)$.  Let $\varphi\in X_\alpha$ be the exact solution of $M\varphi=f$.  Suppose that Assumptions~\ref{ass:logunisolvence}, \ref{ass:continuous}, \ref{ass:discrete}, and \ref{ass:consistency} hold.  Then, for every sufficiently large $n_B$, the exact collocation system \eqref{eq:scheme} has a unique solution $\varphi_{n_B}^H\in V_{n_B,H}$.  Moreover,
\begin{equation*}
        \|\varphi-\varphi_{n_B}^H\|_{X_\beta}
        \le C h_B^{\alpha-\beta}\|\varphi\|_{X_\alpha},
\end{equation*}
where $C$ is independent of $n_B$.
\end{theorem}

\begin{proof}
Assumption~\ref{ass:continuous} guarantees uniqueness and continuous dependence for the exact equation on its range.  The discrete conclusion below uses, in addition, the independent mesh-uniform Assumptions~\ref{ass:discrete} and~\ref{ass:consistency}.  The original and eliminated systems are equivalent because the logarithmic block is nonsingular.  If $v\in V_{n_B,H}$ satisfies the homogeneous collocation equations, then $\widehat E_{n_B}Mv=0$.  Assumption~\ref{ass:discrete} gives $\|v\|_{X_\beta}=0$, hence $v=0$.  Since the number of equations equals $\dim V_{n_B,H}$, the square homogeneous system is nonsingular and the exact collocation system is uniquely solvable for every compatible right-hand side.

Let $w_{n_B}\in V_{n_B,H}$ be the exact-jump approximant from Assumption~\ref{ass:consistency}, and set
\[
        e_{n_B}:=\varphi_{n_B}^H-w_{n_B}\in V_{n_B,H}.
\]
Using $M\varphi=f$ and the collocation equation in the eliminated form, we obtain
\[
        \widehat E_{n_B}Me_{n_B}
        =\widehat E_{n_B}M(\varphi-w_{n_B}).
\]
By scaled discrete stability and scaled consistency,
\[
        \|e_{n_B}\|_{X_\beta}
        \le C_s\|\widehat E_{n_B}M(\varphi-w_{n_B})\|_{\widehat Z_{n_B}^{(\beta)}}
        \le C h_B^{\alpha-\beta}\|\varphi\|_{X_\alpha}.
\]
The triangle inequality and \eqref{eq:scaled_consistency} now give
\[
        \|\varphi-\varphi_{n_B}^H\|_{X_\beta}
        \le \|\varphi-w_{n_B}\|_{X_\beta}+\|e_{n_B}\|_{X_\beta}
        \le C h_B^{\alpha-\beta}\|\varphi\|_{X_\alpha}.
\]
\end{proof}

\begin{remark}
The proof uses the scaled eliminated residual.  This avoids two possible pitfalls: point evaluations of arbitrary logarithmic functions near $\D$ are not required, and the $h_B^{-\beta}$ factor needed to control the $X_\beta$ seminorm of splines is included explicitly in the residual norm.
\end{remark}

\begin{remark}[Scope of Theorem~\ref{thm:main}]\label{rem:conditional_scope}
Theorem~\ref{thm:main} is conditional.  It does not assert that continuous well-posedness automatically implies mesh-uniform stability of the spline collocation matrices, nor that the approximation estimate in $X_\beta$ automatically yields the scaled consistency estimate.  For a particular contour, coefficient set, spline order, phase rule, and family of meshes, Assumptions~\ref{ass:discrete} and~\ref{ass:consistency} require a separate analytical proof.  The finite-dimensional conditioning data reported in Section~\ref{subsec:num_conditioning} provide numerical evidence only for the tested discretization levels and are not a proof of uniform stability as $n_B\to\infty$.
\end{remark}

\section{Quadrature-perturbed collocation}\label{sec:quadrature}

The preceding scheme uses exact values of $S$ and $K$ applied to the trial functions.  In computation these quantities are replaced by quadrature formulas, usually after subtracting the logarithmic or Cauchy singular parts.  We state the perturbation result in the same scaled eliminated norm used in Section~\ref{sec:exactconv}.

Let
\[
        M_{n_B}=M+R_{n_B}
\]
be the operator used in the computed collocation system, where $R_{n_B}$ is the quadrature error operator on $V_{n_B,H}$.  The perturbed scheme is
\begin{equation}\label{eq:perturbedscheme}
        \widehat E_{n_B}M_{n_B}\widetilde\varphi_{n_B}^H
        =\widehat E_{n_B}f .
\end{equation}

\begin{assumption}[Quadrature perturbation]\label{ass:quadrature}
There exists a sequence $\eps_{n_B}\to0$ such that for all $v\in V_{n_B,H}$,
\begin{equation}\label{eq:quadpert}
        \|\widehat E_{n_B}R_{n_B}v\|_{\widehat Z_{n_B}^{(\beta)}}
        \le \eps_{n_B}\|v\|_{X_\beta}.
\end{equation}
\end{assumption}

\begin{theorem}[Perturbed convergence]\label{thm:perturbed}
Under the assumptions of Theorem~\ref{thm:main} and Assumption~\ref{ass:quadrature}, if $n_B$ is sufficiently large so that $C_s\eps_{n_B}<1$, then the perturbed collocation system \eqref{eq:perturbedscheme} has a unique solution $\widetilde\varphi_{n_B}^H\in V_{n_B,H}$.  Moreover, for the approximants $w_{n_B}$ from Assumption~\ref{ass:consistency},
\begin{equation}\label{eq:perturbedestimate}
        \|\varphi-\widetilde\varphi_{n_B}^H\|_{X_\beta}
        \le C\left(
        h_B^{\alpha-\beta}\|\varphi\|_{X_\alpha}
        +\|\widehat E_{n_B}R_{n_B}w_{n_B}\|_{\widehat Z_{n_B}^{(\beta)}}
        \right).
\end{equation}
In particular, if $\eps_{n_B}=O(h_B^{\alpha-\beta})$, then the quadrature-perturbed method has the same asymptotic convergence rate as the exact collocation method.
\end{theorem}

\begin{proof}
For $v\in V_{n_B,H}$, scaled discrete stability gives
\[
        \|v\|_{X_\beta}
        \le C_s\|\widehat E_{n_B}Mv\|_{\widehat Z_{n_B}^{(\beta)}}
        \le C_s\|\widehat E_{n_B}M_{n_B}v\|_{\widehat Z_{n_B}^{(\beta)}}
           +C_s\|\widehat E_{n_B}R_{n_B}v\|_{\widehat Z_{n_B}^{(\beta)}}.
\]
Using \eqref{eq:quadpert},
\[
        (1-C_s\eps_{n_B})\|v\|_{X_\beta}
        \le C_s\|\widehat E_{n_B}M_{n_B}v\|_{\widehat Z_{n_B}^{(\beta)}}.
\]
For $C_s\eps_{n_B}<1$, this is a stability estimate for the perturbed system.  Hence the homogeneous perturbed system has only the zero solution, and the square system is uniquely solvable.

Let
\[
        \widetilde e_{n_B}:=\widetilde\varphi_{n_B}^H-w_{n_B}.
\]
Since $M\varphi=f$, the perturbed equation gives
\[
        \widehat E_{n_B}M_{n_B}\widetilde e_{n_B}
        =\widehat E_{n_B}M(\varphi-w_{n_B})
         -\widehat E_{n_B}R_{n_B}w_{n_B}.
\]
Using the perturbed stability estimate,
\[
        \|\widetilde e_{n_B}\|_{X_\beta}
        \le C\left(
        \|\widehat E_{n_B}M(\varphi-w_{n_B})\|_{\widehat Z_{n_B}^{(\beta)}}
        +\|\widehat E_{n_B}R_{n_B}w_{n_B}\|_{\widehat Z_{n_B}^{(\beta)}}
        \right).
\]
The first term is bounded by the consistency estimate \eqref{eq:scaled_consistency}.  Adding the approximation error gives \eqref{eq:perturbedestimate}.  If \eqref{eq:quadpert} holds with $\eps_{n_B}=O(h_B^{\alpha-\beta})$ and $\|w_{n_B}\|_{X_\beta}\le C\|\varphi\|_{X_\alpha}$, the final assertion follows.
\end{proof}

\begin{remark}[Evaluation of singular integrals]\label{rem:quad_eval}
For trial functions containing relative Heaviside components, the Cauchy integrals contain explicit logarithmic contributions.  A stable implementation should evaluate these logarithmic terms analytically or subtract them before applying a regular quadrature rule.  After this subtraction, the remaining integrands are H\"older continuous on each panel, and composite quadrature rules yield algebraic convergence.  Corrected trapezoidal and hybrid Gauss--trapezoidal quadrature rules for related singular integration problems are discussed in~\cite{KapurRokhlin,Alpert}.  Assumption~\ref{ass:quadrature} is the abstract form of the quadrature accuracy required by the scaled stability argument.
\end{remark}

\section{Numerical experiments}\label{sec:numerical}

This section reports numerical experiments for the regularized B-spline--Heaviside collocation method.  The purpose is fourfold.  First, we verify that the explicit Heaviside enrichment captures the prescribed jumps without visible Gibbs-type oscillations.  Second, we check the decay of the error in untrimmed arcwise discrete analogues of the piecewise H\"older norms used in the analysis.  Third, we document the implementation choices that make the experiments reproducible.  Fourth, we evaluate finite-dimensional diagnostics that are tied directly to the scaled uniform discrete stability and scaled consistency assumptions used in Theorem~\ref{thm:main}.

In both experiments the exact solution is prescribed and the right-hand side is generated by the manufactured-solution identity
\begin{equation}\label{eq:num_rhs_general}
        f=M\varphi=c\varphi+dS\varphi+K\varphi.
\end{equation}
Because the manufactured solutions have nonzero jumps and the coefficient $d$ is not zero at the jump points, the term $S\varphi$ contains logarithmic endpoint contributions.  Thus the right-hand side belongs naturally to the logarithmic target space $Y_\alpha$, not in general to the log-free space $X_\alpha$.  This is precisely the case for which the regularized collocation framework was introduced.

\subsection{Implementation details}\label{subsec:num_implementation}

All computations use the regularized block procedure described in Section~\ref{sec:computational}.  The spline part is represented by periodic cubic B-splines, that is, by splines of order $m=4$.
The discontinuous part is represented by the relative Heaviside basis associated with the prescribed jump set.  The Heaviside amplitudes are determined from logarithmic-coefficient equations, and the point-collocation right-hand side is then regularized by subtracting the explicit contribution of the Heaviside part.  The Cauchy action on B-spline basis functions is evaluated by the singularity-subtracted formula \eqref{eq:comp_S_spline}.  The Cauchy action on relative Heaviside functions is evaluated by the explicit logarithmic arc formula \eqref{eq:comp_S_heaviside}.  The regular operator $K$ is evaluated by panelwise composite trapezoidal quadrature after splitting the parameter interval at the spline knots and at the jump parameters.

The point-collocation nodes coincide with the phase-shifted spline knots,
\[
        \theta_k^B=\rho_{n_B}+(k-1)h_B \pmod{2\pi},
        \qquad h_B=\frac{2\pi}{n_B},
        \qquad k=1,\ldots,n_B.
\]
The phase $\rho_{n_B}\in[0,h_B)$ is selected algorithmically.  The code scans $2000$ candidate phases and chooses the first one satisfying
\[
        \operatorname{dist}(\theta_j^d-\rho_{n_B},h_B\mathbb Z)
        \ge \eta h_B,
        \qquad j=1,\ldots,n_d .
\]
If no candidate satisfies this bound, the phase that maximizes the minimal normalized distance to the jump parameters is used.  In the reported experiments $\eta=0.24$ for Test problem~1 and $\eta=0.15$ for Test problem~2.  This prevents collocation nodes from coinciding with, or becoming too close to, the prescribed jump points.

The quadrature grid is a shifted equidistant midpoint grid,
\[
        \theta_q^Q=\frac{2\pi(q-1)}{N_q}+\frac{\pi}{N_q},
        \qquad q=1,\ldots,N_q.
\]
For the main refinement runs the code uses
\[
        N_q=\max\{N_q^{\min},q_f n_B\}.
\]
In the present runs $N_q^{\min}=131072$, $q_f=240$ for Test problem~1, and $q_f=280$ for Test problem~2.  Hence the ratio $N_q/n_B$ is mesh-dependent until the linear rule $N_q=q_f n_B$ dominates the lower bound $N_q^{\min}$.  The diagnostic grid is independent of the collocation grid and is shifted by half a diagnostic step, so that it does not contain the jump parameters.

In the explicit formula \eqref{eq:comp_S_heaviside}, branch tracking is performed only on the oriented arc that does not contain the evaluation point.  Numerically, for such a nonsingular arc $[\alpha,\beta]$ parametrized by $z(\theta)$, the code samples the path $z(\theta)-t$, unwraps its argument along the oriented path, and evaluates
\[
        \Log(z(\beta)-t)-\Log(z(\alpha)-t)
        =\log|z(\beta)-t|-\log|z(\alpha)-t|
        +i\bigl(\arg_{\rm cont}(z(\beta)-t)-\arg_{\rm cont}(z(\alpha)-t)\bigr).
\]
If $t$ lies on the support arc of $G_j$, the code uses the second line of \eqref{eq:comp_S_heaviside}: it evaluates the nonsingular complementary-arc integral and subtracts it from the full-contour principal value.  Consequently, no logarithmic branch is ever continued through the zero of $z(\theta)-t$.

The manufactured right-hand side is evaluated using the same regularized operator decomposition as the discrete system.  For Test problem~1, the Cauchy transform of the analytic continuous component is evaluated analytically.  For Test problem~2, it is evaluated by the singularity-subtracted quadrature on the high-resolution quadrature grid.  In the quadrature-sensitivity test below, the right-hand side is evaluated on a finer auxiliary grid when the dependence on $N_q$ is examined; this avoids an inverse-crime effect.

The implementation parameters used in the two tests are summarized in Table~\ref{tab:implementation_parameters}.  The diagnostic grid is independent of the collocation grid and does not include the jump parameters themselves.

\begin{table}[H]
\centering
\caption{Implementation parameters used in the numerical experiments.}
\label{tab:implementation_parameters}
\begin{tabular}{ll}
\toprule
Quantity & Choice used in the experiments \\
\midrule
Spline order & $m=4$ cubic periodic B-splines \\
Trial space & $S_{n_B}\oplus \operatorname{span}\{G_2,\ldots,G_{n_d}\}$ \\
Unknowns & $n_B$ spline coefficients and $n_d-1$ Heaviside amplitudes \\
Endpoint equations & logarithmic-coefficient equations \\
Point equations & collocation at nodes separated from $D$ \\
Cauchy term on splines & singularity subtraction, formula \eqref{eq:comp_S_spline} \\
Cauchy term on Heaviside functions & explicit logarithmic arc formula \eqref{eq:comp_S_heaviside} \\
Regular integral & panelwise composite trapezoidal quadrature \\
Linear algebra & solve $R\gamma=b^{\log}$, then $Aa=r$ \\
Precision & complex double precision \\
\bottomrule
\end{tabular}
\end{table}

The point-collocation matrix is dense.  The small logarithmic block has dimension $n_d-1$ and is independent of $n_B$ when the number of jumps is fixed.  The algorithm therefore separates the finite jump-detection part of the computation from the larger spline-collocation part.

\subsection{Error and residual diagnostics}\label{subsec:num_errors}

Let $\Theta_{\rm diag}$ denote the diagnostic set on the continuity arcs and let $t_p=\gamma(\theta_p)$.  The untrimmed maximum error is
\begin{equation*}
        e_\infty^{\rm untr}(n_B)
        =\max_{\theta_p\in\Theta_{\rm diag}}
        |\varphi(t_p)-\varphi_{n_B}^H(t_p)|.
\end{equation*}
The word ``untrimmed'' means that no $O(h_B)$-neighborhood of the jump set is removed.  The jump points themselves are not used as ordinary point-evaluation nodes, because the solution is interpreted there through one-sided values.

To approximate the $X_\beta$ error, the H\"older seminorm is computed arc by arc.  If $\Theta_\ell^{\rm diag}$ is the set of diagnostic parameters on the continuity arc $\Gamma_\ell$, we set
\begin{equation*}
        [e]_{\beta,\Gamma_\ell}^{\rm disc}
        =\max_{\substack{p\ne q\\ \theta_p,\theta_q\in\Theta_\ell^{\rm diag}}}
        \frac{|e(t_p)-e(t_q)|}{|t_p-t_q|^\beta},
        \qquad e=\varphi-\varphi_{n_B}^H.
\end{equation*}
No quotient is formed across a jump.  The reported untrimmed discrete piecewise H\"older-type error is
\begin{equation*}
        e_{X_\beta}^{\rm untr}(n_B)
        =e_\infty^{\rm untr}(n_B)+\max_\ell [e]_{\beta,\Gamma_\ell}^{\rm disc}.
\end{equation*}
When the exact jump amplitudes are known, we also monitor
\begin{equation*}
        e_{\rm jump}(n_B)=\max_{2\le j\le n_d}|\gamma_j-\gamma_j^{\rm exact}|.
\end{equation*}
The implementation additionally computes the point residual and the logarithmic residual,
\begin{equation*}
        \rho_p=\max_{1\le i\le n_B}|(M_{n_B}\varphi_{n_B}^H)(\xi_i)-f(\xi_i)|,
        \qquad
        \rho_{\log}=\max_{2\le j\le n_d}|\Lambda_j(M_{n_B}\varphi_{n_B}^H-f)|.
\end{equation*}
These residuals are algebraic checks of the solved finite system.  They are not substitutes for the error, but they are useful for detecting loss of stability or insufficient quadrature accuracy.

\subsection{Test problem 1: a smooth two-jump mild-logarithmic test}\label{subsec:num_test1}

The first test is posed on the ellipse-like analytic contour
\begin{equation*}
        \Gamma_1=\left\{t=z_1(\theta)=0.75e^{i\theta}+0.25e^{-i\theta}:0\le \theta<2\pi\right\}.
\end{equation*}
The prescribed jump points are
\begin{equation*}
        \theta_1^d=0.70\pi,
        \qquad
        \theta_2^d=1.60\pi,
        \qquad
        t_j^d=z_1(\theta_j^d),\quad j=1,2.
\end{equation*}
The exact solution is
\begin{equation*}
        \varphi(\theta)=
        \begin{cases}
        z_1(\theta)^3+2z_1(\theta), & 0\le \theta\le \theta_1^d\ \text{or}\ \theta_2^d<\theta<2\pi,\\[1mm]
        z_1(\theta)^3+2z_1(\theta)+J, & \theta_1^d<\theta\le \theta_2^d,
        \end{cases}
\end{equation*}
where
\begin{equation*}
        J=0.12+0.04i.
\end{equation*}
Thus $[\varphi]_{t_1^d}=J$ and $[\varphi]_{t_2^d}=-J$.  With the convention $G_2=1$ on the outside arc $[0,\theta_1^d]\cup(\theta_2^d,2\pi)$ and $G_2=0$ on $(\theta_1^d,\theta_2^d]$, one has
\begin{equation*}
        \varphi=\varphi_C+\gamma G_2,
        \qquad
        \varphi_C(\theta)=z_1(\theta)^3+2z_1(\theta)+J,
        \qquad
        \gamma=-J.
\end{equation*}
The continuous component is analytic on the contour.

The coefficients are piecewise H\"older functions with jumps at the same two points:
\begin{equation*}
        c(\theta)=
        \begin{cases}
        2.00+0.06z_1(\theta), & 0\le \theta\le \theta_1^d\ \text{or}\ \theta_2^d<\theta<2\pi,\\[1mm]
        2.15-0.05z_1(\theta), & \theta_1^d<\theta\le \theta_2^d,
        \end{cases}
\end{equation*}
\begin{equation*}
        d(\theta)=0.08
        \begin{cases}
        0.45+0.04z_1(\theta), & 0\le \theta\le \theta_1^d\ \text{or}\ \theta_2^d<\theta<2\pi,\\[1mm]
        0.35-0.03z_1(\theta), & \theta_1^d<\theta\le \theta_2^d.
        \end{cases}
\end{equation*}
The regular kernel is
\begin{equation*}
        h(t,\tau)=0.02\,(t^2+\tau^2).
\end{equation*}
The values $n_B=50,100,200,400$ are used.

\begin{figure}[H]
\centering
\includegraphics[width=0.98\textwidth]{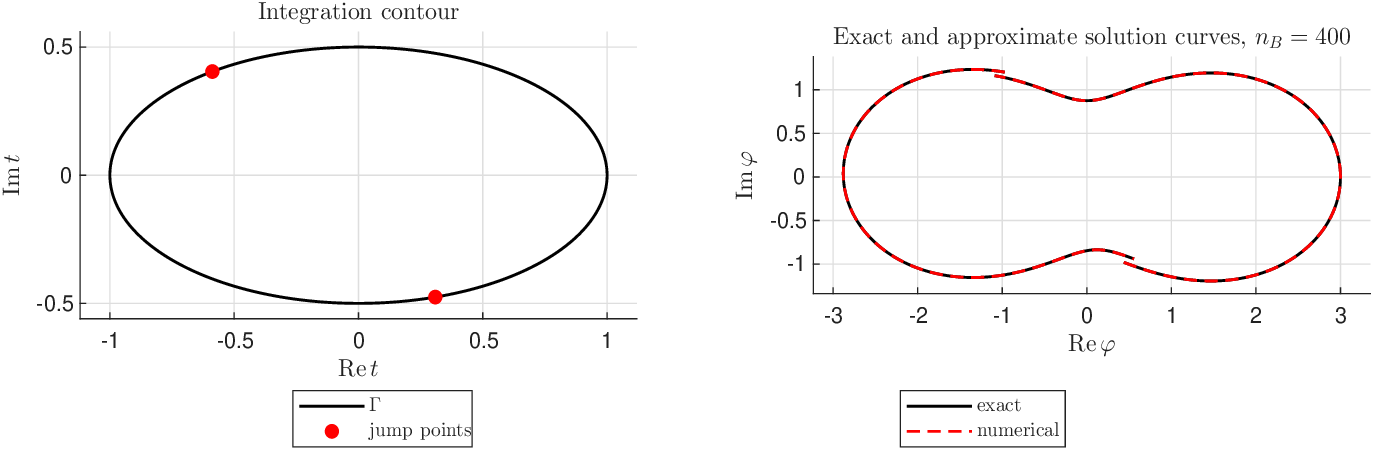}
\caption{Test problem 1: integration contour with the prescribed jump points and exact/approximate solution curves in the complex plane for $n_B=400$.}
\label{fig:num1_contour_solution}
\end{figure}

Figure~\ref{fig:num1_contour_solution} shows that the numerical curve follows the exact solution curve on all smooth branches.  The small breaks correspond to the prescribed jumps.  This representation confirms that the approximation reproduces the geometry of the discontinuous complex-valued solution, not only its real and imaginary parts as scalar functions of the parameter.

\begin{figure}[!htbp]
\centering
\includegraphics[width=0.96\textwidth]{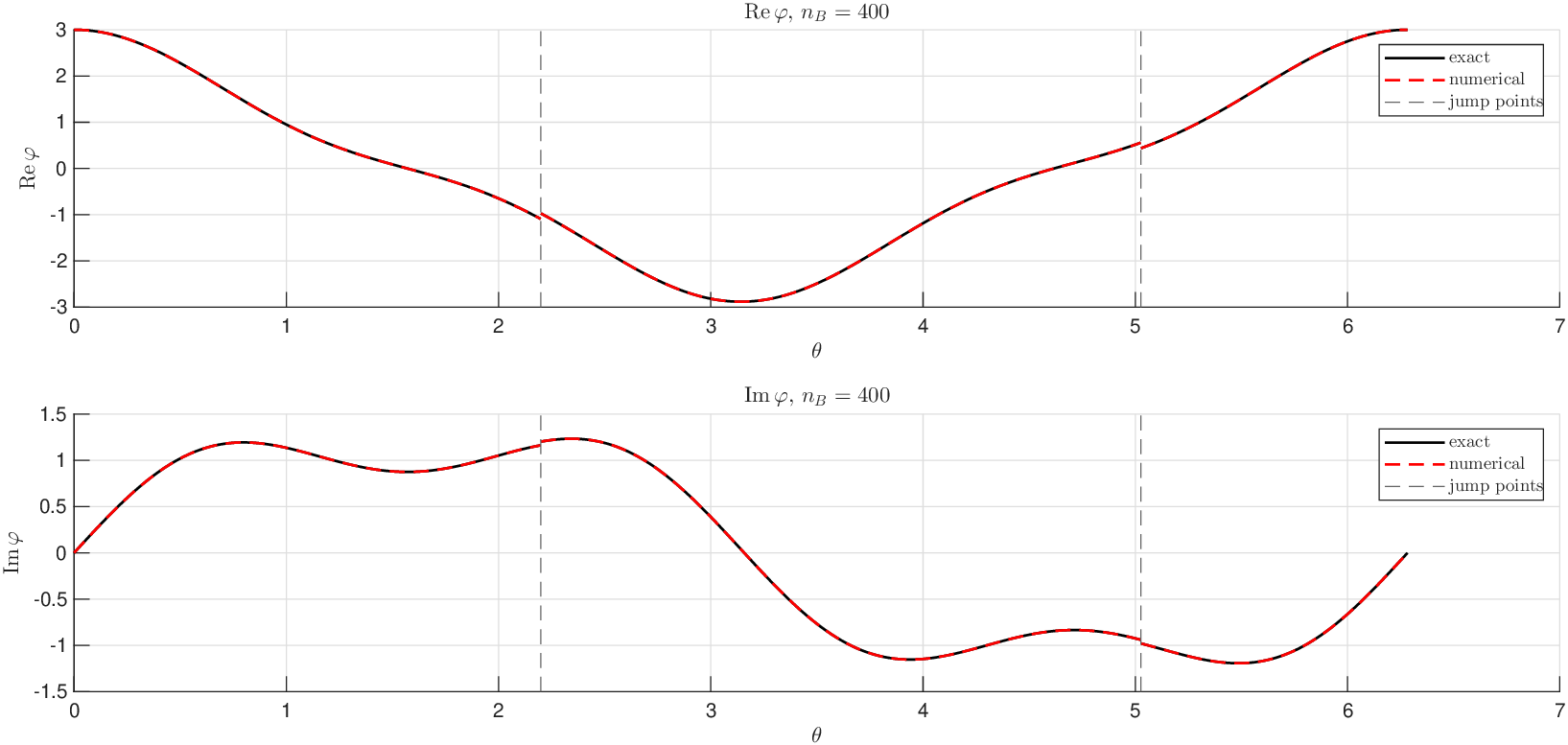}
\caption{Test problem 1: real and imaginary parts of the exact solution and of the B-spline--Heaviside approximation for $n_B=400$.  Dashed vertical lines mark the prescribed jump points.}
\label{fig:num1_solution}
\end{figure}

Figure~\ref{fig:num1_solution} shows that the real and imaginary components are practically indistinguishable on the smooth arcs.  The two prescribed jumps are represented without visible Gibbs-type oscillations.

\begin{figure}[H]
\centering
\includegraphics[width=0.98\textwidth]{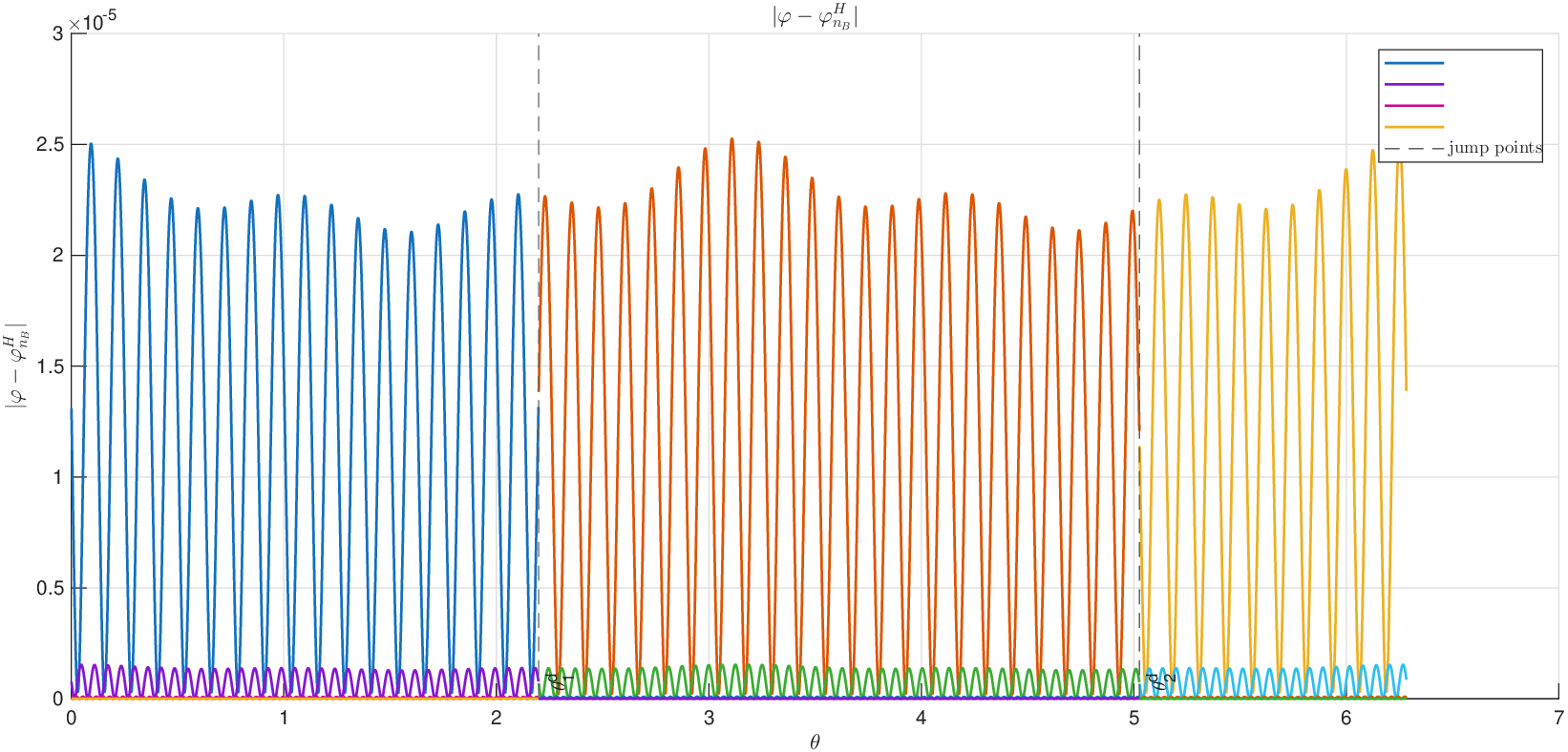}
\caption{Test problem 1: pointwise absolute error $|\varphi-\varphi_{n_B}^H|$ for $n_B=50,100,200,400$.  The vertical axis is scaled by $10^{-5}$.}
\label{fig:num1_pointwise}
\end{figure}

The pointwise error in Figure~\ref{fig:num1_pointwise} is plotted on the smooth arcs, excluding the jump point values themselves.  The visible oscillations are small spline-type oscillations on the smooth arcs, not Gibbs oscillations caused by the discontinuities.  The absence of large spikes near $\theta_1^d$ and $\theta_2^d$ is a direct consequence of the Heaviside enrichment and of the explicit treatment of the logarithmic Cauchy contribution.

\begin{figure}[H]
\centering
\includegraphics[width=0.92\textwidth]{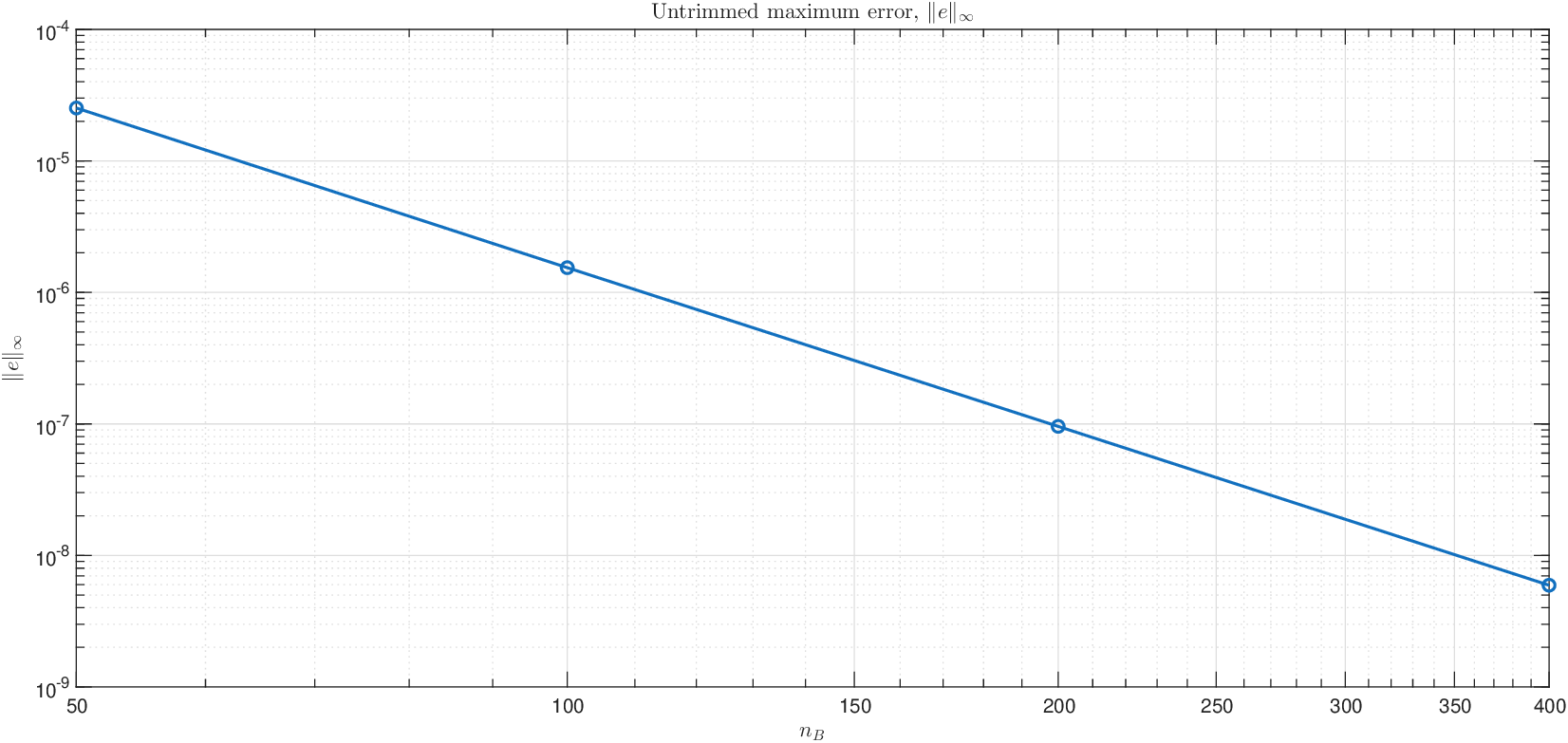}
\caption{Test problem 1: untrimmed maximum error on the full diagnostic grid.  The log-log decay indicates stable algebraic convergence.}
\label{fig:num1_untrimmax}
\end{figure}

The untrimmed maximum errors are reported in Table~\ref{tab:test1_max}.  The observed order is computed by $p=\frac{\log(e(n_B)/e(2n_B))}{\log 2}$.
The values are consistent with the log-log plot in Figure~\ref{fig:num1_untrimmax}.

\begin{table}[H]
\centering
\caption{Test problem 1: untrimmed maximum error and observed order.}
\label{tab:test1_max}
\begin{tabular}{rrrr}
\toprule
$n_B$ & $e_\infty^{\rm untr}$ & reduction factor & observed order \\
\midrule
 50  & $2.5\cdot10^{-5}$ & --   & -- \\
100  & $1.5\cdot10^{-6}$ & $16.7$ & $4.06$ \\
200  & $1.0\cdot10^{-7}$ & $15.0$ & $3.91$ \\
400  & $6.0\cdot10^{-9}$ & $16.7$ & $4.06$ \\
\bottomrule
\end{tabular}
\end{table}

\begin{figure}[H]
\centering
\includegraphics[width=0.92\textwidth]{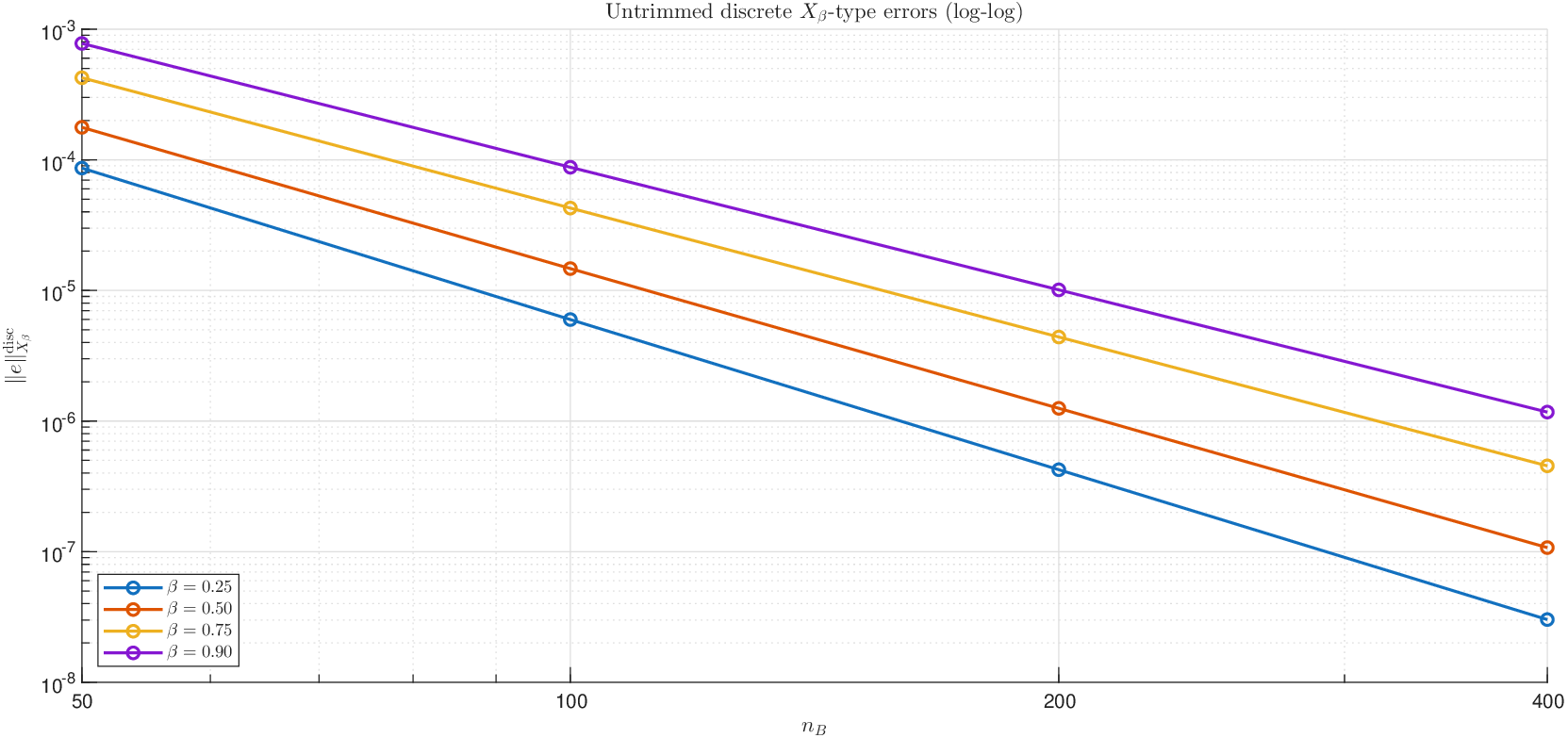}
\caption{Test problem 1: log-log plot of the untrimmed discrete piecewise H\"older-type errors $e_{X_\beta}^{\rm untr}$.  The nearly straight curves indicate algebraic convergence.}
\label{fig:num1_xbeta_untrim_loglog}
\end{figure}

Figure~\ref{fig:num1_xbeta_untrim_loglog} shows the untrimmed discrete piecewise H\"older-type errors for $\beta=0.25,0.50,0.75,0.90$.  The ordering of the curves is natural: larger $\beta$ defines a stronger seminorm and gives a larger error.  The nearly straight curves in the log-log plot indicate algebraic decay.  The observed rates are higher than the minimal theoretical rate $O(h_B^{\alpha-\beta})$ because the continuous component in this test is analytic and cubic splines are used.

\subsection{Test problem 2: a three-jump H\"older-regularity test}\label{subsec:num_test2}

The second test is designed to illustrate a lower-regularity regime in which the behavior is closer to the theoretical dependence on $\alpha-\beta$.  The contour is
\begin{equation*}
        \Gamma_2=\left\{t=z_2(\theta)=e^{i\theta}+0.10e^{5i\theta}:0\le \theta<2\pi\right\}.
\end{equation*}
The prescribed jump parameters are
\begin{equation*}
        \theta_1^d=0.35\pi,
        \qquad
        \theta_2^d=1.05\pi,
        \qquad
        \theta_3^d=1.62\pi,
        \qquad
        t_j^d=z_2(\theta_j^d).
\end{equation*}
Let
\[
        A_1=(\theta_1^d,\theta_2^d],
        \qquad
        A_2=(\theta_2^d,\theta_3^d],
        \qquad
        A_3=(\theta_3^d,2\pi)\cup[0,\theta_1^d].
\]
The continuous component has reference H\"older exponent
$
        \alpha_0=0.72
$
and a H\"older-type singularity at
$
        \theta_c=0.72\pi.
$
It is defined by
\begin{equation*}
        \varphi_C(\theta)
        =0.55z_2(\theta)+0.22z_2(\theta)^2+0.12i\,z_2(\theta)^3
        +(0.18+0.06i)\left|2\sin\frac{\theta-\theta_c}{2}\right|^{\alpha_0}.
\end{equation*}
The exact solution is
\begin{equation*}
        \varphi(\theta)=\varphi_C(\theta)+\gamma_2G_2(\theta)+\gamma_3G_3(\theta),
\end{equation*}
where
\begin{equation*}
        \gamma_2=0.075-0.025i,
        \qquad
        \gamma_3=-0.045+0.055i.
\end{equation*}
Hence
\[
        [\varphi]_{t_1^d}=-\gamma_2-\gamma_3=-0.03-0.03i,
        \qquad
        [\varphi]_{t_2^d}=\gamma_2,
        \qquad
        [\varphi]_{t_3^d}=\gamma_3,
\]
and the total jump is zero.

The coefficients are piecewise H\"older functions with jumps at the same three points:
\begin{equation*}
        c(\theta)=
        \begin{cases}
        1.90+0.045z_2(\theta), & \theta\in A_1,\\[1mm]
        2.06-0.035z_2(\theta)+0.015i, & \theta\in A_2,\\[1mm]
        1.98+0.025i\,z_2(\theta)-0.010i, & \theta\in A_3,
        \end{cases}
\end{equation*}
\begin{equation*}
        d(\theta)=0.055
        \begin{cases}
        0.34+0.020z_2(\theta), & \theta\in A_1,\\[1mm]
        0.29-0.018z_2(\theta), & \theta\in A_2,\\[1mm]
        0.32+0.012i\,z_2(\theta), & \theta\in A_3.
        \end{cases}
\end{equation*}
The regular kernel is
\begin{equation*}
        h(t,\tau)=0.015\,(t^2+\tau^2).
\end{equation*}
The right-hand side is generated by \eqref{eq:num_rhs_general}.  The values $n_B=80,160,320,640$ are used, and the tested H\"older exponents are
\[
        \beta=0.20,
        \qquad
        \beta=0.35,
        \qquad
        \beta=0.50,
        \qquad
        \beta=0.65,
\]
all smaller than $\alpha_0=0.72$.

\begin{figure}[!htbp]
\centering
\includegraphics[width=0.96\textwidth]{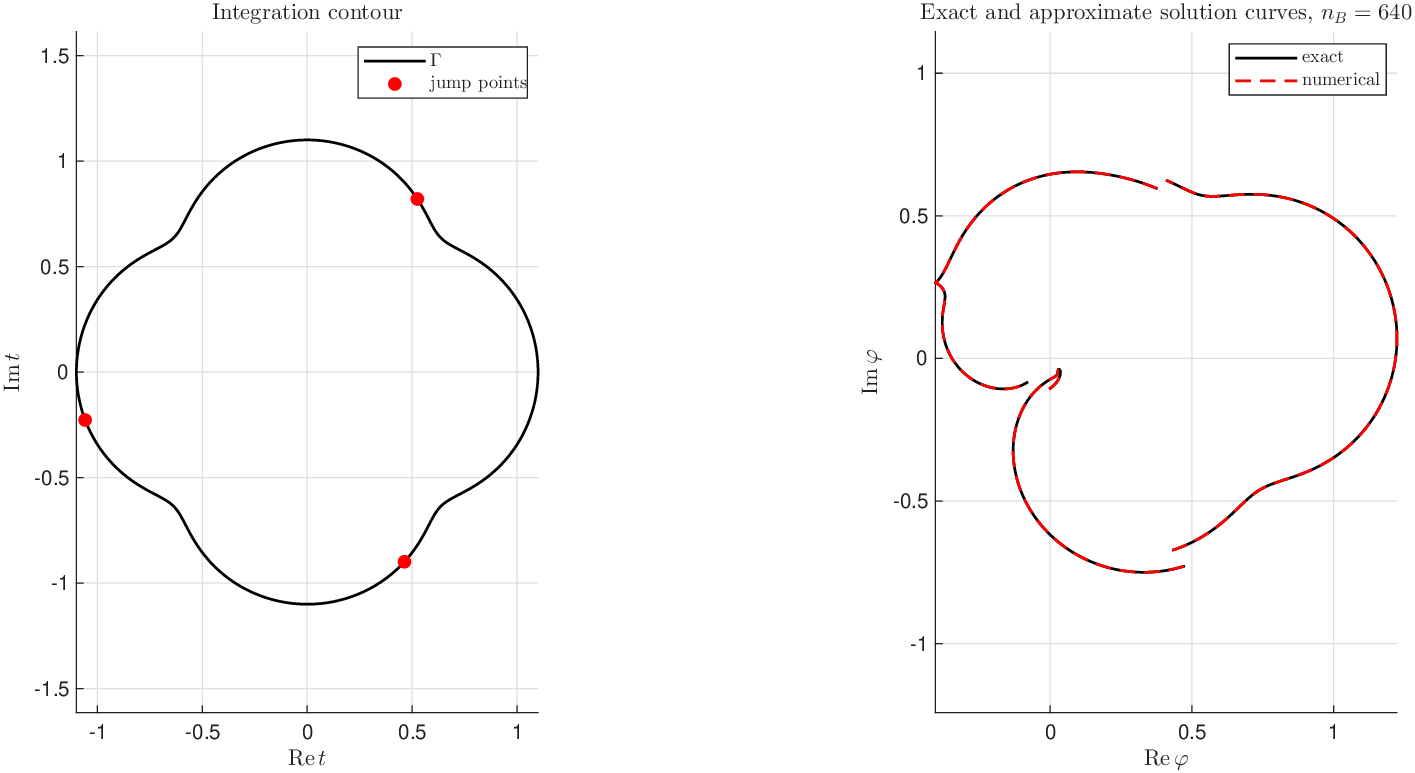}
\caption{Test problem 2: the integration contour $\Gamma_2$ with the three jump points, and the exact and approximate solution curves in the $(\operatorname{Re}\varphi,\operatorname{Im}\varphi)$ plane for $n_B=640$.}
\label{fig:num2_contour_solution}
\end{figure}

Figure~\ref{fig:num2_contour_solution} shows that the exact and numerical solution curves nearly coincide.  The small breaks correspond to the prescribed jumps.  This confirms that the jump structure is captured by the relative Heaviside enrichment without visible spurious loops.

\begin{figure}[H]
\centering
\includegraphics[width=0.98\textwidth]{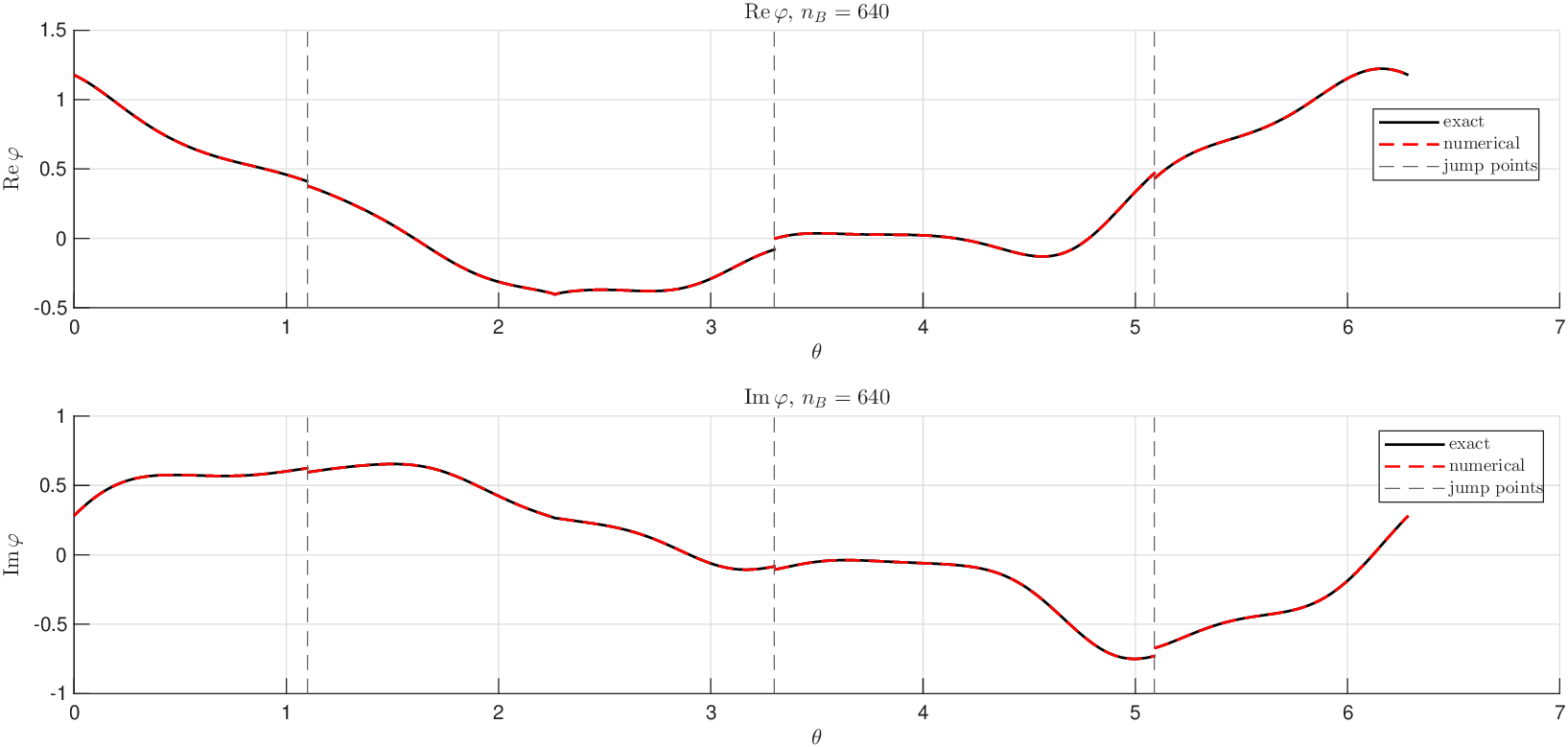}
\caption{Test problem 2: real and imaginary parts of the exact and approximate solutions for $n_B=640$.  Dashed vertical lines mark the three jump parameters.}
\label{fig:num2_solution_theta}
\end{figure}

Figure~\ref{fig:num2_solution_theta} shows that the numerical solution follows the exact one closely on all continuity arcs and reproduces the correct one-sided branches at the jump points.  No visible instability is observed near the discontinuities.

\begin{figure}[!htbp]
\centering
\includegraphics[width=0.88\textwidth]{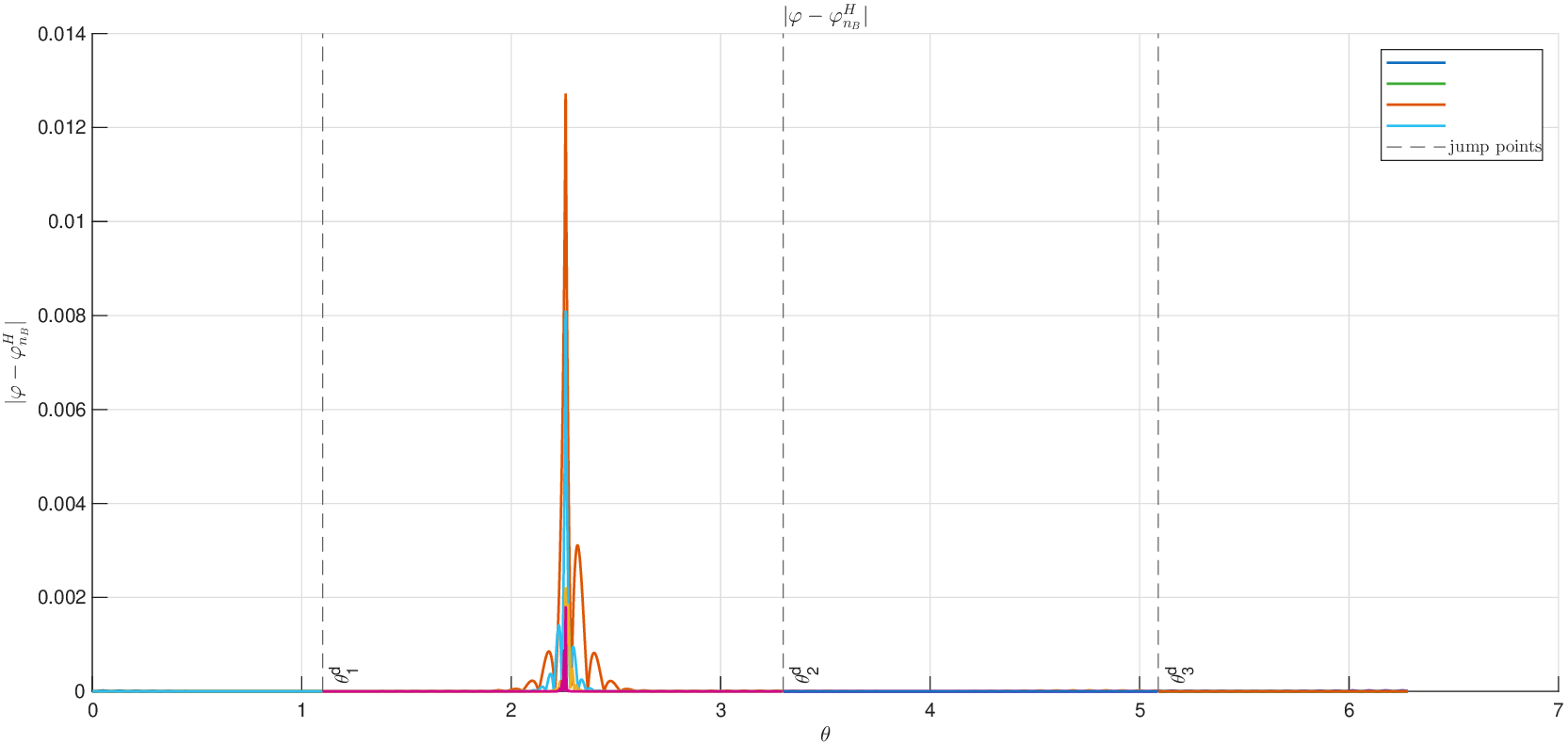}
\caption{Test problem 2: pointwise absolute error for $n_B=80,160,320,640$.  The dominant peak is localized near the H\"older singularity of the continuous component, not at the prescribed jump points.}
\label{fig:num2_pointwise}
\end{figure}

The pointwise error in Figure~\ref{fig:num2_pointwise} is dominated by the H\"older singularity at $\theta_c=0.72\pi$, not by the prescribed jump points.  This is the intended behavior: the Heaviside enrichment removes the jump-induced approximation difficulty, while the remaining error is governed by the regularity of the continuous component.

The untrimmed maximum errors are summarized in Table~\ref{tab:test2_max}.  The decay is slower than in Test problem~1, as expected for a continuous component with limited H\"older regularity.

\begin{table}[H]
\centering
\caption{Test problem 2: untrimmed maximum error and observed order.}
\label{tab:test2_max}
\begin{tabular}{rrrr}
\toprule
$n_B$ & $e_\infty^{\rm untr}$ & reduction factor & observed order \\
\midrule
 80  & $1.25\cdot10^{-2}$ & --   & -- \\
160  & $8.0\cdot10^{-3}$  & $1.56$ & $0.64$ \\
320  & $2.3\cdot10^{-3}$  & $3.48$ & $1.80$ \\
640  & $1.7\cdot10^{-3}$  & $1.35$ & $0.44$ \\
\bottomrule
\end{tabular}
\end{table}

\begin{figure}[!htbp]
\centering
\includegraphics[width=0.88\textwidth]{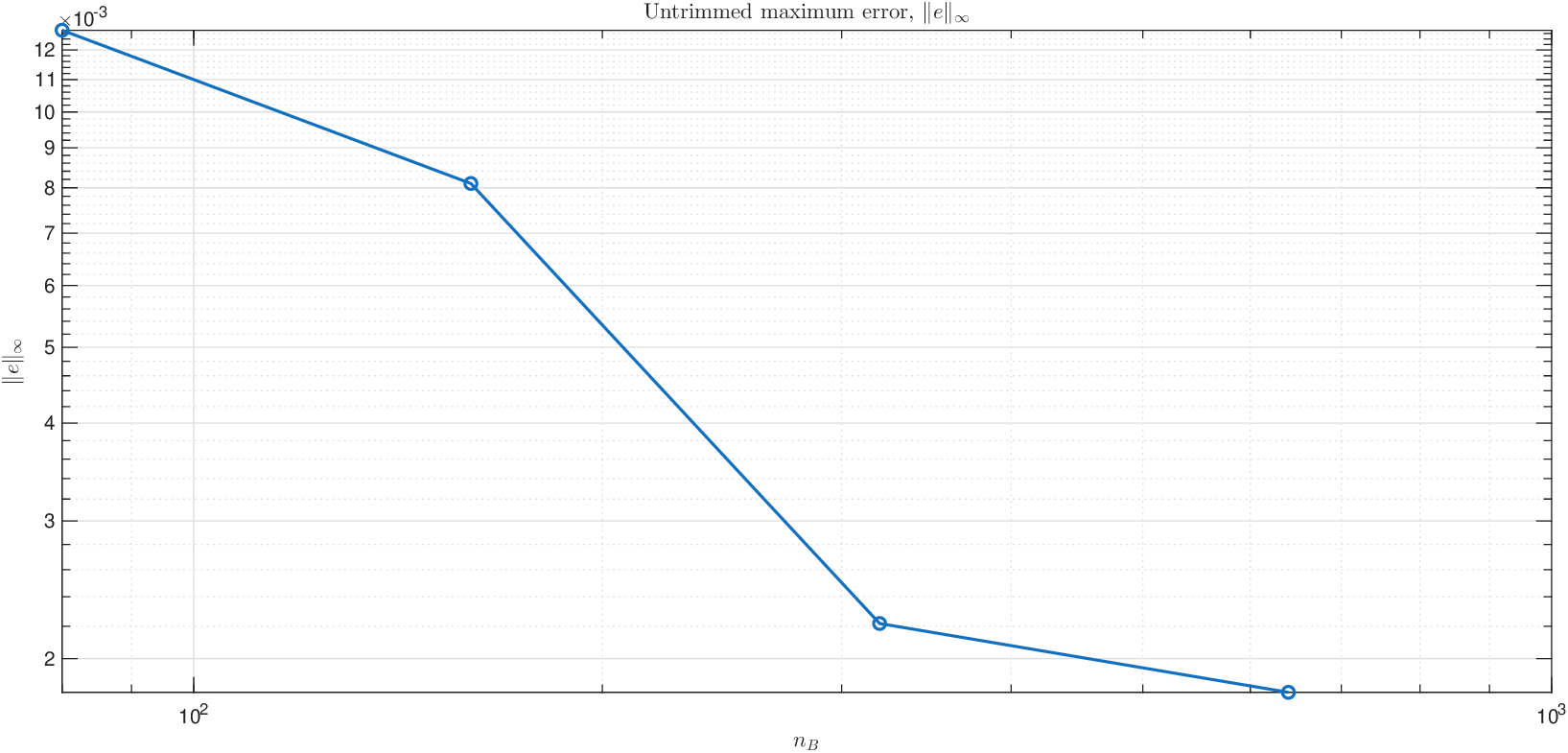}
\caption{Test problem 2: untrimmed maximum error on the full diagnostic grid.  The slower decay reflects the H\"older singularity in the continuous component.}
\label{fig:num2_untrimmax}
\end{figure}

\begin{figure}[!htbp]
\centering
\includegraphics[width=0.88\textwidth]{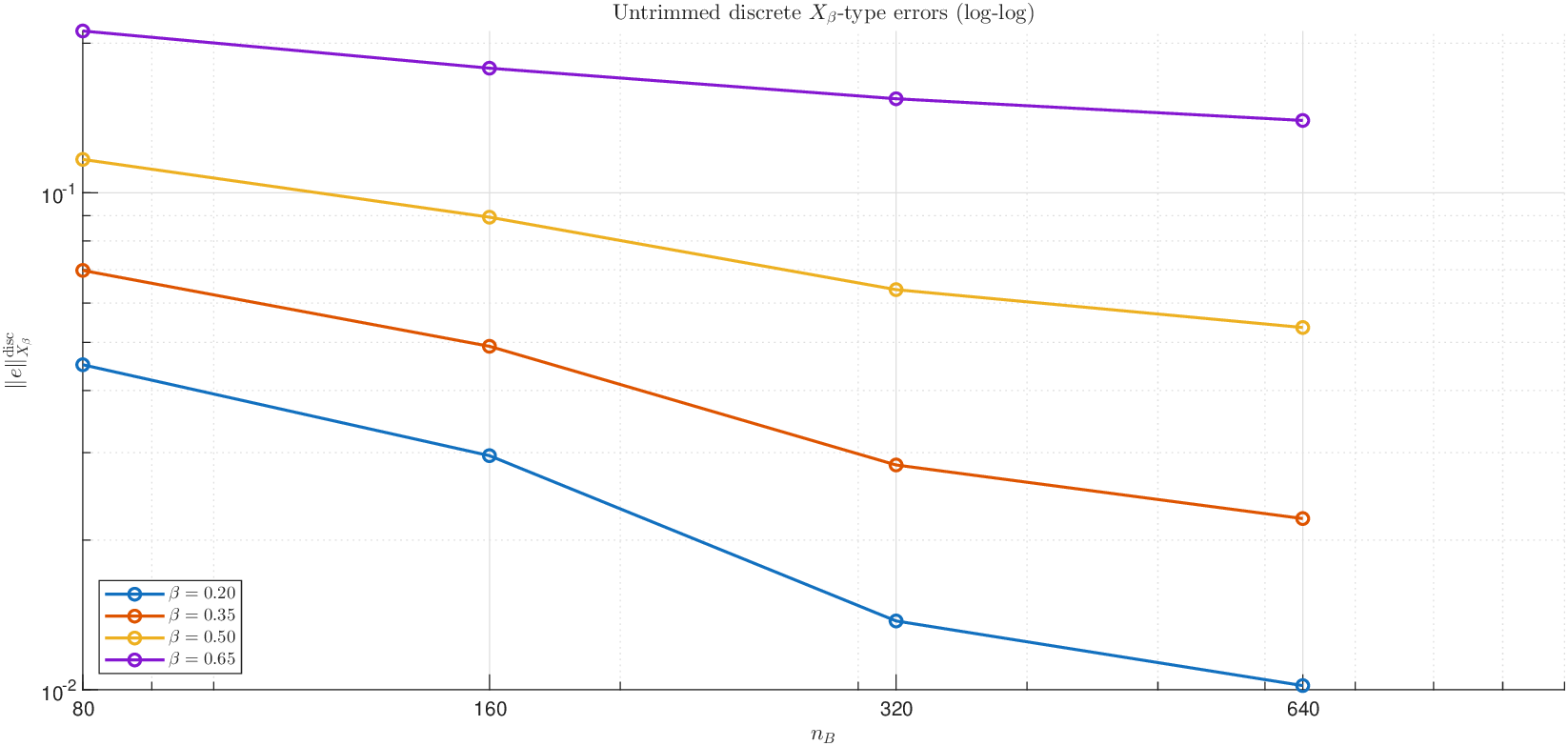}
\caption{Test problem 2: log-log plot of the untrimmed discrete piecewise H\"older-type errors.  The algebraic decay is slower for larger $\beta$, in agreement with the dependence on $\alpha-\beta$.}
\label{fig:num2_xbeta_untrim_loglog}
\end{figure}

The maximum-error plot in Figure~\ref{fig:num2_untrimmax} shows the same slower convergence behaviour as Table~\ref{tab:test2_max}.  The error is controlled by the H\"older cusp of the continuous component rather than by the prescribed jumps.  The untrimmed discrete piecewise H\"older-type errors are shown in Figure~\ref{fig:num2_xbeta_untrim_loglog}.  Their ordering
\[
        e_{X_{0.20}}^{\rm untr}
        <e_{X_{0.35}}^{\rm untr}
        <e_{X_{0.50}}^{\rm untr}
        <e_{X_{0.65}}^{\rm untr}
\]
is consistent with the fact that larger $\beta$ gives a stronger seminorm.  The log-log plot indicates algebraic convergence, and the decay becomes slower as $\beta$ approaches $\alpha_0=0.72$.  This agrees with the theoretical mechanism in Theorem~\ref{thm:main}, where the approximation order is governed by $\alpha-\beta$.

\FloatBarrier

To obtain a more robust rate estimate for the lower-regularity test, we also ran an extended refinement study for $\beta=0.50$ with
\[
        n_B=80,160,320,640,1280,2560.
\]
For the larger systems, sparse B-spline quadrature matrices and target batching were used.  The results are reported in Table~\ref{tab:extended_test2}.  Since the continuous component contains a localized H\"older cusp, consecutive two-grid rates are not fully regular in the pre-asymptotic range.  We therefore estimate the rate by a least-squares regression of $\log e_{X_\beta}^{\rm untr}$ against $\log n_B$ on the last three and on the last four successful refinement levels.

\begin{table}[H]
\centering
\caption{Extended refinement study for Test problem~2 with $\beta=0.50$.}
\label{tab:extended_test2}
\begin{tabular}{rrrrr}
\toprule
$n_B$ & $N_q$ & $e_\infty^{\rm untr}$ & $e_{X_\beta}^{\rm untr}$ & $\rho_p$ \\
\midrule
  80 &  32768 & $1.27202\cdot10^{-2}$ & $1.16754\cdot10^{-1}$ & $2.245\cdot10^{-15}$ \\
 160 &  32768 & $8.10125\cdot10^{-3}$ & $8.93318\cdot10^{-2}$ & $2.121\cdot10^{-15}$ \\
 320 &  32768 & $2.21840\cdot10^{-3}$ & $6.38534\cdot10^{-2}$ & $3.998\cdot10^{-15}$ \\
 640 &  51200 & $1.81150\cdot10^{-3}$ & $5.35769\cdot10^{-2}$ & $6.707\cdot10^{-15}$ \\
1280 & 102400 & $1.13176\cdot10^{-3}$ & $3.99431\cdot10^{-2}$ & $1.272\cdot10^{-14}$ \\
2560 & 204800 & $1.03393\cdot10^{-3}$ & $4.13034\cdot10^{-2}$ & $1.776\cdot10^{-14}$ \\
\bottomrule
\end{tabular}

\vspace{1mm}
For all rows, $\rho_{\log}=2.168\cdot10^{-19}$.
\end{table}

For this test, the reference H\"older exponent of the continuous component is $\alpha_0=0.72$, so the theoretical rate in the $X_\beta$ estimate is governed by
\[
        \alpha_0-\beta=0.22.
\]
The regression rate computed from $e_{X_\beta}^{\rm untr}$ on the last three refinement levels is $0.188$, while the rate computed on the last four refinement levels is $0.231$.  These values are close to the reference value $0.22$ and support the predicted dependence on $\alpha-\beta$.  The irregular consecutive rates are consistent with a pre-asymptotic regime caused by the localized H\"older cusp.

\subsection{Direct comparison with classical B-spline collocation}\label{subsec:num_comparison}

To isolate the effect of the Heaviside enrichment, we also solved the same point-collocation equations using only the continuous periodic B-spline space $S_{n_B}$.  In this comparison, the logarithmic endpoint equations are not used and the trial function has zero jumps.  Consequently, the method can only approximate the discontinuous solution indirectly by a continuous spline.  Table~\ref{tab:direct_comparison} compares the proposed method with the classical B-spline approximation for Test problem~2 with $n_B=640$ and $\beta=0.50$.

\begin{table}[H]
\centering
\caption{Direct comparison between the proposed B-spline--Heaviside method and classical B-spline collocation for Test problem~2, with $n_B=640$ and $\beta=0.50$.}
\label{tab:direct_comparison}
\begin{tabular}{lccc}
\toprule
Method & $e_\infty^{\rm untr}$ & $e_{X_\beta}^{\rm untr}$ & $e_{\rm jump}$ \\
\midrule
B-spline--Heaviside & $1.81150\cdot10^{-3}$ & $5.35769\cdot10^{-2}$ & $1.388\cdot10^{-17}$ \\
Classical B-spline  & $6.60258\cdot10^{-2}$ & $6.88666\cdot10^{-1}$ & $7.906\cdot10^{-2}$ \\
\bottomrule
\end{tabular}
\end{table}

\begin{figure}[!htbp]
\centering
\includegraphics[width=0.88\textwidth]{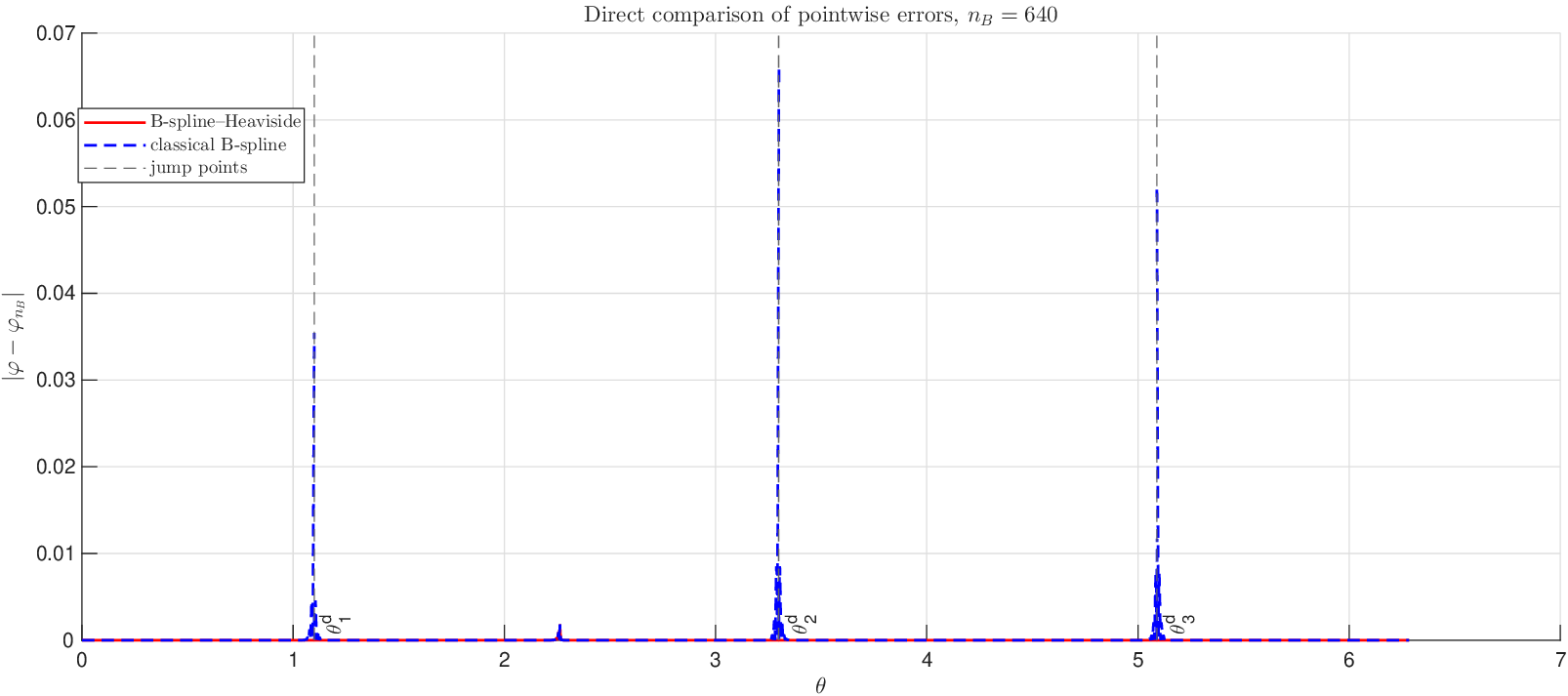}
\caption{Direct comparison between the proposed B-spline--Heaviside method and classical B-spline collocation for Test problem~2 with $n_B=640$. The classical B-spline approximation develops localized pointwise error spikes near the prescribed jump points, while the enriched method suppresses these spikes by representing the jumps explicitly.}
\label{fig:direct_comparison_pointwise_errors}
\end{figure}

The numerical advantage of representing the jumps explicitly is substantial.  The pointwise-error comparison in Figure~\ref{fig:direct_comparison_pointwise_errors} shows that, although the classical B-spline approximation may look visually close to the exact solution away from the discontinuity points, its pointwise error develops pronounced spikes near the prescribed jumps.  The proposed B-spline--Heaviside method removes these jump-induced spikes and recovers the Heaviside amplitudes to machine precision.  In this test, the untrimmed maximum error of the classical B-spline approximation is about $36$ times larger than that of the B-spline--Heaviside method, while the discrete $X_\beta$-type error is about $13$ times larger.  This confirms that the enrichment is not merely a cosmetic modification of the approximation space, but is essential for reproducing the piecewise H\"older structure of the solution.

\FloatBarrier

\subsection{Quadrature sensitivity}\label{subsec:num_quadrature}

The convergence theorem for the perturbed scheme requires the quadrature contribution to be smaller than, or at least of the same order as, the discretization error.  We tested this by fixing $n_B=320$ in Test problem~2, taking $\beta=0.50$, increasing the ratio $N_q/n_B$, and monitoring
\[
        e_\infty^{\rm untr},
        \qquad
        e_{X_\beta}^{\rm untr},
        \qquad
        \rho_p,
        \qquad
        \rho_{\log}.
\]
The right-hand side in this test was evaluated on a finer auxiliary quadrature grid.  This separates the quadrature error in the assembled system from the error used to manufacture the data.

\begin{table}[H]
\centering
\caption{Quadrature-sensitivity diagnostics for Test problem~2 with $n_B=320$ and $\beta=0.50$.}
\label{tab:quadrature_sensitivity}
\begin{tabular}{rrrrr}
\toprule
$N_q/n_B$ & $N_q$ & $e_\infty^{\rm untr}$ & $e_{X_\beta}^{\rm untr}$ & $\rho_p$ \\
\midrule
 20.0 &   6400 & $2.21839\cdot 10^{-3}$ & $6.38534\cdot 10^{-2}$ & $4.003\cdot 10^{-15}$ \\
 40.0 &  12800 & $2.21839\cdot 10^{-3}$ & $6.38534\cdot 10^{-2}$ & $6.233\cdot 10^{-15}$ \\
 80.0 &  25600 & $2.21839\cdot 10^{-3}$ & $6.38534\cdot 10^{-2}$ & $4.003\cdot 10^{-15}$ \\
160.0 &  51200 & $2.21839\cdot 10^{-3}$ & $6.38534\cdot 10^{-2}$ & $4.453\cdot 10^{-15}$ \\
280.0 &  89600 & $2.21839\cdot 10^{-3}$ & $6.38534\cdot 10^{-2}$ & $4.011\cdot 10^{-15}$ \\
\bottomrule
\end{tabular}

\vspace{1mm}
For all rows, $\rho_{\log}=2.168\cdot10^{-19}$.
\end{table}

Table~\ref{tab:quadrature_sensitivity} shows that the reported errors are essentially unchanged when the quadrature refinement is increased from $N_q/n_B=20$ to $N_q/n_B=280$.  The point residual remains at the level of roundoff error, and the logarithmic residual is constant at approximately $2.168\cdot10^{-19}$.  Hence, in the tested range, the quadrature error is well below the discretization error and does not determine the observed convergence behavior.

\subsection{Conditioning and computational cost}\label{subsec:num_conditioning}

The regularized method leads to two principal algebraic objects: the small logarithmic block $R$ and the dense point-collocation matrix $A$.  We therefore monitor
\begin{equation*}
        \kappa(R),
        \qquad
        \kappa(A),
        \qquad
        \kappa(\mathcal A),
        \qquad
        \mathcal A=
        \begin{pmatrix}
        A&H\\0&R
        \end{pmatrix}.
\end{equation*}
The condition numbers reported below are computed after the same row and column scaling used in the linear solve.  Thus they describe the conditioning of the actually solved algebraic systems rather than an arbitrarily scaled matrix representation.  The assembly time includes the construction of $R,A,H,r$ and of the right-hand side, while the solve time includes the solution of $R\gamma=b^{\log}$ and $Aa=r$.

\begin{table}[H]
\centering
\caption{Conditioning and computational cost diagnostics for Test problem~2.  Condition numbers are computed after the same row and column scaling used in the linear solve.}
\label{tab:conditioning_cost}
\begin{tabular}{rrrrrrrr}
\toprule
$n_B$ & $N_q/n_B$ & $\kappa(R)$ & $\kappa(A)$ & $\kappa(\mathcal A)$ & assembly (s) & solve (s) & $e_{\rm jump}$ \\
\midrule
 80  & 1638.4 & $1.07039$ & $3.02984$ & $5.07094\cdot10^{1}$ & 0.6092 & 0.0500 & $1.388\cdot10^{-17}$ \\
160  &  819.2 & $1.07039$ & $3.03158$ & $9.81793\cdot10^{1}$ & 0.7958 & 0.0066 & $1.388\cdot10^{-17}$ \\
320  &  409.6 & $1.07039$ & $3.03252$ & $1.93976\cdot10^{2}$ & 1.7094 & 0.0134 & $1.388\cdot10^{-17}$ \\
640  &  280.0 & $1.07039$ & $3.03292$ & $3.83933\cdot10^{2}$ & 6.0960 & 0.0527 & $1.388\cdot10^{-17}$ \\
\bottomrule
\end{tabular}

\vspace{1mm}
For all rows, $\rho_{\log}=2.168\cdot10^{-19}$, while $\rho_p$ ranges between $1.555\cdot10^{-15}$ and $1.110\cdot10^{-14}$.
\end{table}

Table~\ref{tab:conditioning_cost} indicates good algebraic conditioning of the two principal blocks in the tested range.  The logarithmic block $R$ has condition number close to one, and the row/column-scaled point-collocation matrix $A$ has condition number close to three for all tested values of $n_B$.  The condition number of the full block matrix $\mathcal A$ increases with $n_B$, but this growth does not lead to a visible loss of accuracy: both the point residual and the logarithmic residual remain close to roundoff level, and the jump amplitudes are recovered to essentially machine precision.

Condition numbers alone do not verify Assumption~\ref{ass:discrete}.  In particular, multiplication of a matrix by a scalar leaves its condition number unchanged, whereas the scaled residual norm in \eqref{eq:scaled_residual_norm} contains the factor $h_B^{-\beta}$.  Moreover, a small residual of the solved algebraic system checks the linear solve but not the independent consistency estimate required in Assumption~\ref{ass:consistency}.  For this reason we next report diagnostics constructed directly from the inverse maps and from exact-jump approximants.

\subsection{Scaled stability and consistency diagnostics}\label{subsec:num_scaled_diagnostics}

For the logarithmic block we compute
\begin{equation*}
        C_R(n_B):=\|R^{-1}\|_{\ell^\infty\to\ell^\infty}.
\end{equation*}
This is the finite-dimensional constant entering condition (S1) of Proposition~\ref{prop:block_scaled}.  Its magnitude is problem-dependent and is not required to be close to one; the relevant issue is whether it remains bounded under refinement.

Let $B_{\rm diag}$ be the matrix that evaluates the spline basis on the independent diagnostic grid.  Since $A$ maps spline coefficients to the regularized point residual, we also compute
\begin{equation*}
        C_A^{\rm num}(n_B):=\|B_{\rm diag}A^{-1}\|_{\ell^\infty\to\ell^\infty}.
\end{equation*}
This quantity is a discrete-grid proxy for the maximum-norm stability constant in condition (S2): it measures the amplification from point residual values to spline values on the diagnostic grid.  It is more directly related to (S2) than $\kappa(A)$, although it is still a finite-grid numerical indicator rather than a proof of a continuous supremum-norm bound.

To test scaled consistency independently of the solved collocation solution, we construct
\begin{equation*}
        w_{n_B}=I_{n_B}\varphi_C+\sum_{j=2}^{n_d}\gamma_j^{\rm exact}G_j,
\end{equation*}
where $I_{n_B}\varphi_C$ is the periodic cubic-spline interpolant of the continuous component at the phase-shifted knots.  Thus $w_{n_B}$ has exactly the same jumps as $\varphi$.  We record
\begin{align*}
        r_{\rm cons}^{\rm point}(n_B)
        &:=\max_{1\le k\le n_B}
          |M(\varphi-w_{n_B})(t_k^B)|,\\
        \eta_{\rm cons}(n_B)
        &:=h_B^{-\alpha}r_{\rm cons}^{\rm point}(n_B),\\
        r_{\rm cons}^{\log}(n_B)
        &:=\max_{2\le j\le n_d}
          |\Lambda_j(M(\varphi-w_{n_B}))|.
\end{align*}
Boundedness of $\eta_{\rm cons}$ is the numerical signature of the point-residual requirement $r_{\rm cons}^{\rm point}=O(h_B^\alpha)$ in Assumption~\ref{ass:consistency}.  Since the jumps are inserted exactly, $r_{\rm cons}^{\log}$ should vanish up to roundoff.

\begin{table}[H]
\centering
\small
\caption{Scaled stability and consistency diagnostics for Test problem~1 ($\alpha=0.99$).}
\label{tab:scaled_diag_test1}
\begin{tabular}{rrrrrr}
\toprule
$n_B$ & $C_R$ & $C_A^{\rm num}$ & $r_{\rm cons}^{\rm point}$ & $\eta_{\rm cons}$ & $r_{\rm cons}^{\log}$ \\
\midrule
 50  & $2.701232\cdot10^{1}$ & $7.907312\cdot10^{-1}$ & $6.380512\cdot10^{-7}$  & $4.973221\cdot10^{-6}$ & $0$ \\
100  & $2.701232\cdot10^{1}$ & $7.965901\cdot10^{-1}$ & $3.576208\cdot10^{-8}$  & $5.536364\cdot10^{-7}$ & $0$ \\
200  & $2.701232\cdot10^{1}$ & $8.025286\cdot10^{-1}$ & $2.120934\cdot10^{-9}$  & $6.521518\cdot10^{-8}$ & $0$ \\
400  & $2.701232\cdot10^{1}$ & $8.064914\cdot10^{-1}$ & $1.291905\cdot10^{-10}$ & $7.889905\cdot10^{-9}$ & $0$ \\
\bottomrule
\end{tabular}
\end{table}

For Test problem~1, $C_R(n_B)$ is independent of the mesh and $C_A^{\rm num}(n_B)$ varies only from approximately $0.791$ to $0.806$.  Hence the two direct stability proxies remain uniformly bounded over the reported refinements.  The raw point consistency residual decreases from $6.38\cdot10^{-7}$ to $1.29\cdot10^{-10}$, and the normalized indicator $\eta_{\rm cons}$ decreases by almost three orders of magnitude.  This is stronger than the boundedness required by the $O(h_B^{0.99})$ consistency condition and reflects the analytic continuous component in this test.  The logarithmic consistency residual is zero to the displayed precision.

\begin{table}[H]
\centering
\small
\caption{Scaled stability and consistency diagnostics for Test problem~2 ($\alpha=0.72$).}
\label{tab:scaled_diag_test2}
\begin{tabular}{rrrrrr}
\toprule
$n_B$ & $C_R$ & $C_A^{\rm num}$ & $r_{\rm cons}^{\rm point}$ & $\eta_{\rm cons}$ & $r_{\rm cons}^{\log}$ \\
\midrule
 80  & $5.882594\cdot10^{1}$ & $8.455160\cdot10^{-1}$ & $1.230974\cdot10^{-4}$ & $7.687482\cdot10^{-4}$ & $0$ \\
160  & $5.882594\cdot10^{1}$ & $8.490706\cdot10^{-1}$ & $6.441444\cdot10^{-5}$ & $6.626133\cdot10^{-4}$ & $0$ \\
320  & $5.882594\cdot10^{1}$ & $8.511080\cdot10^{-1}$ & $4.243852\cdot10^{-5}$ & $7.190824\cdot10^{-4}$ & $0$ \\
640  & $5.882594\cdot10^{1}$ & $8.541291\cdot10^{-1}$ & $2.812269\cdot10^{-5}$ & $7.849046\cdot10^{-4}$ & $0$ \\
\bottomrule
\end{tabular}
\end{table}

For Test problem~2, $C_R(n_B)=58.82594$ is again mesh-independent, while $C_A^{\rm num}(n_B)$ remains in the narrow interval $[0.8455,0.8542]$.  The point consistency residual decreases under refinement, and the scaled indicator remains between $6.63\cdot10^{-4}$ and $7.85\cdot10^{-4}$.  Its bounded, mildly nonmonotone behavior is consistent with
\[
        r_{\rm cons}^{\rm point}(n_B)=O(h_B^{0.72}),
\]
which is the rate required by Assumption~\ref{ass:consistency} for this rough H\"older test.  The small variation of $\eta_{\rm cons}$ is compatible with the pre-asymptotic fluctuations already observed in the $X_\beta$ errors.  Again, the logarithmic consistency residual vanishes to the displayed precision because the exact jumps are built into $w_{n_B}$.

Taken together, Tables~\ref{tab:scaled_diag_test1} and~\ref{tab:scaled_diag_test2} give numerical evidence for the two implementation-dependent assumptions at the tested discretization levels: the inverse-norm proxies associated with (S1)--(S2) remain bounded, and the exact-jump interpolants satisfy the expected scaled point-consistency behavior.  These tables do not prove uniform bounds for all $n_B$, but they test the quantities that occur directly in the theoretical estimates, rather than relying only on condition numbers or on residuals of the already solved system.

\subsection{Summary of the numerical evidence}\label{subsec:num_summary}

The two tests illustrate complementary aspects of the method.  Test problem~1 has an analytic continuous component and shows high practical convergence rates.  Test problem~2 contains a H\"older singularity in the continuous component and therefore exhibits slower algebraic decay, closer to the theoretical dependence on $\alpha-\beta$.  In both tests the prescribed jumps are represented explicitly and no visible Gibbs-type oscillations appear at the jump points.  The untrimmed arcwise diagnostics are closer to the theoretical $PH^\beta(\Gamma,D)$ norm than trimmed diagnostics: they use the full diagnostic grid on each continuity arc and avoid only the jump point values themselves.

The direct scaled diagnostics complement these error measurements.  In both problems $\|R^{-1}\|_\infty$ is independent of $n_B$, the numerical spline stability proxy $\|B_{\rm diag}A^{-1}\|_\infty$ remains nearly constant, and the exact-jump logarithmic consistency residual vanishes.  For the rough test the normalized point-consistency indicator remains bounded, whereas for the smooth test it decreases rapidly.  Thus the experiments support, at the tested mesh levels, both the approximation behavior and the implementation-dependent stability and consistency hypotheses used by the conditional convergence theorem.  They remain numerical evidence rather than an analytical proof of uniformity as $n_B\to\infty$.

\section{Conclusions}\label{sec:conclusion}

We have developed a B-spline--Heaviside collocation framework for Cauchy singular integral equations with piecewise H\"older solutions on a closed contour.  The solution is measured in the space $X_\alpha=PH^\alpha(\Ga,\D)$,
where jumps are finite and explicitly represented.  The residual and the right-hand side are measured in the logarithmically enlarged space
$Y_\alpha=PH^\alpha_{\log,*}(\Ga,\D)$,
which contains the logarithmic singularities produced by the Cauchy operator acting on discontinuous functions.

The use of lateral H\"older logarithmic coefficients is essential.  It makes $Y_\alpha$ a module over $X_\alpha$, so that $dS\varphi\in Y_\alpha$ whenever $d\in X_\alpha$ and $\varphi\in X_\alpha$.  The right-hand side is therefore naturally allowed to belong to $Y_\alpha$.  The smaller assumption $f\in X_\alpha$ is possible only under compatibility conditions that eliminate the logarithmic coefficients, for example when the solution has no jump at points where $d$ is nonzero.

On a closed contour, the discontinuous enrichment must be formulated with relative Heaviside functions.  These functions have zero total jump and form a basis of the finite-dimensional jump space.  The enriched approximation space
\[
        V_{n_B,H}=S_{n_B}\oplus\spanop\{G_2,\ldots,G_{n_d}\}
\]
has dimension $n_B+n_d-1$ and approximates every $\varphi\in X_\alpha$ in the weaker norm $X_\beta$, $0<\beta<\alpha<1$, with rate
\[
        \inf_{v\in V_{n_B,H}}
        \|\varphi-v\|_{X_\beta}
        \le C h_B^{\alpha-\beta}\|\varphi\|_{X_\alpha}.
\]

The collocation scheme is formulated through bounded regularized functionals: point evaluations at mesh nodes separated from the jump set and logarithmic-coefficient functionals at the jump points.  After eliminating the Heaviside contribution determined by the logarithmic equations, the point residual is measured in the scaled norm dictated by the spline inverse inequality.  Under continuous stability, scaled uniform discrete stability, and scaled consistency of exact-jump approximants, the exact collocation solution exists uniquely for all sufficiently fine meshes and satisfies the same $O(h_B^{\alpha-\beta})$ convergence rate.  An abstract perturbation theorem shows that quadrature-based implementations preserve the rate when the quadrature error is sufficiently small in the scaled eliminated residual norm.

The numerical evidence shows that, for the tested regularized matrices, the method behaves stably and achieves the predicted arcwise convergence.  In particular, the logarithmic block is well conditioned, the Heaviside amplitudes are recovered to machine precision, and the comparison with classical B-spline collocation confirms that explicit jump representation removes the localized error spikes produced near discontinuity points.  More directly, the computed quantities $\|R^{-1}\|_\infty$ and $\|B_{\rm diag}A^{-1}\|_\infty$ remain bounded over all reported refinements.  The normalized exact-jump consistency indicator $h_B^{-\alpha}\max_k|M(\varphi-w_{n_B})(t_k^B)|$ is bounded for the rough H\"older test and decreases rapidly for the smooth test, while the corresponding logarithmic residual is zero to the displayed precision.

The remaining implementation-dependent issue is the analytical verification of the zero-jump spline stability and the scaled consistency estimate for specific choices of spline order, phase-shifted meshes, logarithmic-coefficient functionals, and quadrature rules.  These properties are not consequences of the continuous mapping theory and are not proved here for arbitrary coefficient sets.  The present formulation isolates this question, separates it from the logarithmic singularities, and provides a self-contained functional framework in which such a verification can be carried out for a concrete B-spline--Heaviside collocation family.

\end{document}